\documentstyle[12pt,amscd]{amsart}
\newtheorem{thm}[equation]{Theorem}
\numberwithin{equation}{section}
\newtheorem{cor}[equation]{Corollary}
\newtheorem{expl}[equation]{Example}
\newtheorem{rmk}[equation]{Remark}

\newtheorem{lem}[equation]{Lemma}

\newtheorem{diag}[equation]{Diagram}
\newtheorem{prop}[equation]{Proposition}
\newtheorem{tab}[equation]{Table}

\newtheorem{stmt}[equation]{Statement}

\begin{document}
\raggedbottom \voffset=-.7truein \hoffset=0truein \vsize=8truein \hsize=6truein \textheight=8truein \textwidth=6truein \baselineskip=18truept
\def\mapright#1{\ \smash{\mathop{\longrightarrow}\limits^{#1}}\ }
\def\mapleft#1{\smash{\mathop{\longleftarrow}\limits^{#1}}}
\def\mapup#1{\Big\uparrow\rlap{$\vcenter {\hbox {$#1$}}$}}
\def\mapdown#1{\Big\downarrow\rlap{$\vcenter {\hbox {$\ssize{#1}$}}$}}
\def\mapne#1{\nearrow\rlap{$\vcenter {\hbox {$#1$}}$}}
\def\mapse#1{\searrow\rlap{$\vcenter {\hbox {$\ssize{#1}$}}$}}
\def\mapr#1{\smash{\mathop{\rightarrow}\limits^{#1}}}
\def\l{\lambda}
\def\TK{T_{K/2}}
\def\ss{\smallskip}
\def\vp{v_1^{-1}\pi}
\def\at{{\widetilde\alpha}}
\def\sm{\wedge}
\def\la{\langle}
\def\ra{\rangle}
\def\on{\operatorname}
\def\spin{\on{Spin}}
\def\kbar{{\overline k}}
\def\lbar{{\overline{\ell}}}
\def\qed{\quad\rule{8pt}{8pt}\bigskip}
\def\ssize{\scriptstyle}
\def\a{\alpha}
\def\Z{{\Bbb Z}}
\def\F{{\Bbb F}}
\def\im{\on{im}}
\def\ct{\widetilde{C}}
\def\gdt{\widetilde{\on{gd}}}
\def\ext{\on{Ext}}
\def\sq{\on{Sq}}
\def\eps{\epsilon}
\def\ar#1{\stackrel {#1}{\rightarrow}}
\def\br{{\bold R}}
\def\bC{{\bold C}}
\def\bh{{\bold H}}
\def\bx{{\bold x}}
\def\si{\sigma}
\def\Ebar{{\overline E}}
\def\dbar{{\overline d}}
\def\Sum{\sum}
\def\tfrac{\textstyle\frac}
\def\tb{\textstyle\binom}
\def\Si{\Sigma}
\def\w{\wedge}
\def\equ{\begin{equation}}
\def\b{\beta}
\def\G{\Gamma}
\def\g{\gamma}
\def\psit{\widetilde{\Psi}}
\def\tht{\widetilde{\Theta}}
\def\psiu{{\underline{\Psi}}}
\def\thu{{\underline{\Theta}}}
\def\aee{A_{\text{ee}}}
\def\aeo{A_{\text{eo}}}
\def\aoo{A_{\text{oo}}}
\def\aoe{A_{\text{oe}}}
\def\jbar{{\bold{j}}}
\def\Qh{\hat\Q_7}
\def\endeq{\end{equation}}
\def\sn{S^{2n+1}}
\def\zp{\bold Z_p}
\def\A{{\cal A}}
\def\P{{\cal P}}
\def\cj{{\cal J}}
\def\zt{{\bold Z}_2}
\def\bs{{\bold s}}
\def\bof{{\bold f}}
\def\bq{{\Bbb Q}}
\def\bQ{{\Bbb Q}}
\def\be{{\bold e}}
\def\Hom{\on{Hom}}
\def\ker{\on{ker}}
\def\coker{\on{coker}}
\def\da{\downarrow}
\def\colim{\operatornamewithlimits{colim}}
\def\zphat{\bz_2^\wedge}
\def\io{\iota}
\def\Om{\Omega}
\def\u{{\cal U}}
\def\e{{\cal E}}
\def\exp{\on{exp}}
\def\wbar{{\overline w}}
\def\xbar{{\overline x}}
\def\ybar{{\overline y}}
\def\zbar{{\overline z}}
\def\ebar{{\bold {e}}}
\def\nbar{{\overline n}}
\def\rbar{{\overline r}}
\def\et{{\widetilde E}}
\def\ni{\noindent}
\def\coef{\on{coef}}
\def\den{\on{den}}
\def\lcm{\on{l.c.m.}}
\def\vi{v_1^{-1}}
\def\ot{\otimes}
\def\psibar{{\overline\psi}}
\def\mhat{{\hat m}}
\def\exc{\on{exc}}
\def\ms{\medskip}
\def\ehat{{\hat e}}
\def\etao{{\eta_{\text{od}}}}
\def\etae{{\eta_{\text{ev}}}}
\def\dirlim{\operatornamewithlimits{dirlim}}
\def\gt{\widetilde{L}}
\def\lt{\widetilde{\lambda}}
\def\st{\widetilde{s}}
\def\ft{\widetilde{f}}
\def\Ft{\widetilde{F}}
\def\el{\ell}
\def\om{\omega}
\def\sgd{\on{sgd}}
\def\lfl{\lfloor}
\def\rfl{\rfloor}
\def\ord{\on{ord}}
\def\gd{{\on{gd}}}
\def\rk{{{\on{rk}}_2}}
\def\nbar{{\overline{n}}}
\def\lg{{\on{lg}}}
\def\N{{\Bbb N}}
\def\Z{{\Bbb Z}}
\def\Q{{\Bbb Q}}
\def\R{{\Bbb R}}
\def\C{{\Bbb C}}
\def\mo{\on{mod}}
\def\Zh{{\hat\Z_7}}
\def\vexp{v_1^{-1}\exp}
\def\Remark{\noindent{\it  Remark}}
\title[$p$-compact groups]
{Homotopy type and $v_1$-periodic homotopy groups of $p$-compact groups}

\author{Donald M. Davis}
\address{Lehigh University\\Bethlehem, PA 18015}
\email{dmd1@@lehigh.edu}
\thanks{The author is grateful to Kasper Andersen, John Harper, Hirokazu Nishinobu, and Clarence Wilkerson for valuable suggestions
related to this project, each of which is specifically described in the paper.}
\date{October 25, 2007}

\maketitle
\begin{abstract} We determine the $v_1$-periodic homotopy groups of all irreducible $p$-compact groups $(BX,X)$. In the most
difficult, modular, cases, we follow a direct path
from their associated invariant polynomials to these homotopy groups. We show that, if $p$ is odd,
every irreducible $p$-compact group has $X$ of the homotopy type of a product of explicit  spaces related to
 $p$-completed Lie groups.
\end{abstract}

\section{Introduction}\label{intro} In \cite{4auth} and \cite{Gro}, the classification of
irreducible $p$-compact groups was completed. This family of spaces extends the family of ($p$-completions of) compact simple Lie groups.
The $v_1$-periodic homotopy groups of any space $X$, denoted $\vp_*(X)_{(p)}$, are a localization of the portion of the homotopy groups
detected by $K$-theory; they were defined in \cite{DM}. In
\cite{Reprth} and \cite{DHHA}, the author completed the determination of the $v_1$-periodic homotopy groups of all compact simple Lie groups.
Here we do the same
 for all the remaining irreducible $p$-compact groups.\footnote{If the groups $\vp_i(X)$ are finite, then
$p$-completion induces an isomorphism $\vp_*(X)\to\vp_*(X_p)$. (\cite[p.1252]{Bo})}

Recall that a $p$-compact group (\cite{DW}) is a pair $(BX,X)$ such that $BX$ is $p$-complete and $X=\Omega BX$ with $H^*(X;\F_p)$ finite. Thus $BX$
determines $X$ and contains more structure than does $X$. The homotopy type and homotopy groups of $X$ do not take into account this extra
structure nor the multiplication on $X$. We show that, if $p$ is odd,
every irreducible $p$-compact group has $X$ of the homotopy type of a product of explicit  spaces related to
 $p$-completed Lie groups.

According to \cite[1.1,11.1]{4auth} and \cite{Gro}, the irreducible $p$-compact groups correspond  to
 compact simple Lie groups\footnote{Cases in which distinct
compact Lie groups give rise to equivalent $p$-compact groups are discussed in \cite[11.4]{4auth}.} and the $p$-adic reflection groups listed in \cite[Table 1]{And} for which the character field is
strictly larger than $\Q$.
See \cite[pp.430-431]{CE} and \cite[p.165]{Kane} for other listings of reflection groups.  We use the usual notation $((BX_n)_p,(X_n)_p)$, where $n$
is the Shephard-Todd numbering (\cite{ST} or any of the previously-mentioned tables) and $p$ is the prime associated to the completion.

We will divide our discussion into four families of cases:
\begin{enumerate}
\item The compact simple Lie groups---infinite family 1, part of infinite family 2, and cases 28, 35, 36, 37 in the Shephard-Todd list.
\item The rest of the infinite families numbered 2a, 2b, and 3. \item The nonmodular special cases, in which $p$ does not divide the order of the reflection group. This is cases 4-27
and 29-34. \item The modular cases, in which $p$ divides the order
of the reflection group. These are  cases $(X_{12})_3$, $(X_{24})_2$, $(X_{29})_5$, $(X_{31})_5$, and $(X_{34})_7$. (Actually, we include
$(X_{12})_3$ in Case (3) along with the nonmodular cases, and the Dwyer-Wilkerson space $(X_{24})_2$ was handled in \cite{DI4}.)
\end{enumerate}

Here is a brief summary of what we accomplish in each case. The author feels that his contributions here are nil in case (1)\footnote{But he accomplished
much  in these cases in earlier papers such as \cite{DHHA}, \cite{Reprth}, and \cite{DSU}.}, minuscule in case (2), modest
in case (3), and significant in case (4).
\begin{enumerate}
\item Spaces $X_1$, $X_{28}$, $X_{35}$, $X_{36}$, and $X_{37}$ are, respectively, $SU(n)$, $F_4$, $E_6$, $E_7$, and $E_8$. These are $p$-compact
groups for all primes $p$, although for small primes they were excluded by Clark and Ewing (\cite{CE}) because $H^*(BX;\F_p)$ is not a polynomial algebra.
The exceptional Lie group $G_2$ is the case $m=6$ in infinite family 2b. The spaces $SO(n)$, Spin$(n)$, and $Sp(n)$ appear in the infinite family
2a with $m=2$.
Simplification of the homotopy types of many of these, when $p$ is odd, to products of spheres and spherically-resolved spaces was obtained in
\cite[(8.1),8.1]{MNT}. The $v_1$-periodic homotopy groups of these spaces were computed in \cite{DSU}, \cite{BDSU}, \cite{Reprth}, \cite{DHHA},
 \cite{BDMi}, and other papers. We will say no more about these cases. \item In Section \ref{infsec}, we use work of Castellana and Broto-Moller to show that the spaces in
the infinite families can be decomposed, up to homotopy, as products of factors of $p$-completions of unitary groups, spheres, and sphere bundles over
spheres. See \ref{m|p-1}, \ref{rmk1}, and \ref{m|p+1} for the specific results. \item In Table \ref{gtbl}, we list the homotopy types of all
cases $(X_n)_p$ which are not products of spheres. There are 31 such cases. In each case, we give the homotopy type as a product of spheres and
spaces which are spherically resolved with $\a_1$ attaching maps. In Remark \ref{rmk2}, we discuss the easily-computed $v_1$-periodic homotopy
groups of these spaces. \item The most novel part of the paper is the determination of the $v_1$-periodic homotopy groups of $(X_{29})_5$,
 $(X_{31})_5$, and $(X_{34})_7$. We introduce a direct, but nontrivial, path from the invariant polynomials to the $v_1$-periodic homotopy groups.
En route, we determine the Adams operations in $K^*(BX;\hat\Z_p)$ and $K^*(X;\hat\Z_p)$. In the case of $(X_{34})_7$, we give new explicit formulas for
the invariant polynomials. We prove in
\ref{29conj} (resp. \ref{34conj}) that the homotopy type of $(X_{29})_5$ (resp. $(X_{34})_7$) is directly related to $SU(20)$ (resp. $SU(42))$.
We prove in \ref{X31E8} that $(X_{31})_5$ has the homotopy type of the 5-completion of a factor of $E_8$. These latter homotopy-type results
and their proofs were pointed out and explained to the author by John Harper.
\end{enumerate}

\section{Infinite families 2 and 3}\label{infsec}
Family 3 consists of $p$-completed\footnote{All of our spaces are completed at an appropriate prime $p$. This will not always be present in our
notation. For example, we will often write $SU(n)$ when we really mean its $p$-completion.} spheres $S^{2m-1}$ with $p\equiv1$ mod $m$, which is
a loop space due to work of Sullivan (\cite{Sul}). The groups $\vp_*(S^{2m-1})_{(p)}$, originally due to Mahowald ($p=2$) and Thompson ($p$
odd), are given in \cite[4.2]{Dsur}.

Family 2 consists of spaces $X(m,r,n)$ where $m>1$, $r|m$, and $n>1$. The ``degrees"  of $X(m,r,n)$
are $m,2m,\ldots,(n-1)m,\frac mrn$.  These are the degrees
of invariant polynomials under a group action used in defining the space. The Clark-Ewing table doubles the degrees to form the ``type," as
these doubled degrees are the degrees of generators of $H^*(BX;\F_p)$ in the cases which they consider.
For most\footnote{According to \cite{Not}, the only exclusions are certain compact Lie groups when $p$ is very small.} of the
irreducible $p$-compact groups $X$, $H^*(X;\F_p)$ is an exterior
algebra on classes of grading $2d-1$, where $d$ ranges over the degrees. Family 2b consists of spaces $X(m,r,n)$ in which $n=2$ and $r=m$, while
family 2a is all other cases. The reason that these are separated is that 2b has more applicable primes. Indeed, for family 2a, there are
$p$-compact groups when $p\equiv1$ mod $m$, while for family 2b, these exist when $p\equiv\pm1$ mod $m$, and also $p=2$ if $m=4$ or 6, and $p=3$
if $m=3$ or 6. The case $m=6$ in family 2b is the exceptional Lie group $G_2$. Note that all primes work when $m=6$. The case $(p=2,m=4)$ has
$X=Sp(2)$, while $(p=3,m=3)$ has $X=SU(3)$ or $PSU(3)$, the projective unitary group. In this case, there are two inequivalent $p$-compact groups
corresponding to the same $Q_p$-reflection group; however, since $SU(3)\to PSU(3)$ is a 3-fold covering space, they have isomorphic
$v_1$-periodic homotopy groups.

The following results of Broto and Moller (\cite{BM}) and Castellana (\cite{Cas}) will be useful. They deal with the homotopy fixed-point space
$X^{hG}$ when $G$ acts on a space of the same homotopy type as a space $X$. Here and throughout, $C_m$ denotes a cyclic group of order $m$, and
$U(N)$ is the $p$-completion of a unitary group.
\begin{thm}\label{BM1} $(\cite[5.2,5.12]{BM})$ If $m|(p-1)$, $0\le s<m$, and $n>0$, then
$$U(mn+s)^{hC_m}\simeq X(m,1,n)$$ and is a factor in a product decomposition of $U(mn+s)$.\end{thm}
\begin{thm}\label{BM2} $(\cite[5.2,5.14]{BM})$ If $m|(p-1)$,
$m\ge2$, $r>1$, and $n\ge2$,  then $$X(m,r,n)^{hC_m}\simeq X(m,1,n-1)$$ and is a factor in a product decomposition of $X(m,r,n)$.\end{thm}

\begin{cor}\label{m|p-1} If $m|(p-1)$ and $r>1$, then
$$X(m,r,n)\simeq X(m,1,n-1)\times S^{2n\frac mr-1}$$
and $X(m,1,n-1)$ is a factor in a product decomposition of $U(m(n-1))$.\end{cor}

Here $X(m,1,1)$ is interpreted as $S^{2m-1}$.

\begin{pf} We use Theorem \ref{BM2} to get the first factor. By the Kunneth Theorem,
the other factor must have the same $\F_p$-cohomology as $S^{2n\frac mr-1}$, and hence must have the same homotopy type as this sphere. Now we
apply Theorem \ref{BM1} to complete the proof.\end{pf}

\begin{rmk}{\rm Our Corollary \ref{m|p-1} appears as \cite[1.4]{Cas}, except that she has an apparent typo regarding the dimension of the sphere.
Also, neither she nor \cite{BM} have the restriction $r>1$, but it seems that the result is false for $r=1$, since by induction it would imply
that $X(m,1,n)$ is a product of spheres, which is not usually true.}\end{rmk}

\begin{rmk}{\rm Let $p$ be odd. By \cite{MNT}, for any $N$, $p$-completed $SU(N)$ splits as a product of $(p-1)$ spaces, each of which
has $H^*(-;\F_p)$ an exterior algebra on odd dimensional classes of dimensions $b$, $b+q,\ldots,b+tq$, for some integers $b$ and $t$.
Here and throughout $q=2(p-1)$. Our space $X(m,1,n-1)$ will be a product of
$(p-1)/m$ of these spaces for $SU(m(n-1))$. The $v_1$-periodic homotopy groups of these spaces can be read off from those of $SU(m(n-1))$, since
the $(p-1)$ factors have $v_1$-periodic homotopy groups in nonoverlapping dimensions. Thus, to the extent that \cite{DSU} is viewed as being a
satisfactory description of $\vp_*(SU(n))_{(p)}$,\footnote{\cite[1.4]{DSU} states that $\vp_{2k}(SU(n))_{(p)}$ is a cyclic $p$-group with
exponent $\min(\nu_p(j!S(k,j)):j\ge n)$, where $S(-,-)$ denotes the Stirling number of the second kind. In \cite{DY}, more tractable formulas
were obtained if $n\le p^2-p+1$. Here and throughout, $\nu_p(-)$ is the exponent of $p$.} Corollary \ref{m|p-1} gives $\vp_*(X(m,r,n))_{(p)}$
provided $m|(p-1)$.}\label{rmk1}\end{rmk}

\begin{expl} Let $p=7$. Then $X(2,2,6)\simeq X(2,1,5)\times S^{11}$. There is a product decomposition
$$(SU(10))_7\simeq B(3,15)\times B(5,17)\times B(7,19)\times S^9\times S^{11}\times S^{13},$$
where $B(2n+1,2n+13)$ denotes a 7-completed $S^{2n+1}$-bundle over $S^{2n+13}$ with attaching map $\a_1$. Then
$$X(2,1,5)\simeq B(3,15)\times B(7,19)\times S^{11}.$$
\end{expl}

What remains for Family 2 is the cases 2b when $m|(p+1)$. These are the spaces $X(m,m,2)$ with $m|(p+1)$.  Let $B(3,2p+1)$ denote the $p$-completion of
an $S^3$-bundle
over $S^{2p+1}$ with attaching map $\a_1$.
\begin{thm}\label{m|p+1} If $m|(p+1)$, then
$$X(m,m,2)\simeq\begin{cases}B(3,2p+1)&m=p+1\\
S^3\times S^{2m-1}&m<p+1.\end{cases}$$
\end{thm}
\begin{pf} Let $X=X(m,m,2)$ with $m|(p+1)$. Then $H^*(X;\F_p)=\Lambda[x_3,x_{2m-1}]$.
If $m<p+1$, then by the unstable Adams spectral sequence (\cite{6auth}), both classes $x_3$ and $x_{2m-1}$ are spherical. Indeed, the $E_2$-term
begins with towers in dimensions 3 and $2m-1$ emanating from filtration 0, and no possible differentials. See Diagram \ref{sphch}.
Because $X$ is
an $H$-space, the maps $S^3\to X$ and $S^{2m-1}\to X$ yield a map $S^3\times S^{2m-1}\to X$, and it is a $p$-equivalence by Whitehead's Theorem.

On the other hand, suppose $m=p+1$. We will show that $\P^1(x_3)=x_{2p+1}$. It then follows from \cite[7.1]{MNT} that there is a $p$-equivalence
$B(3,2p+1)\to X$.

To see that $\P^1(x_3)=x_{2p+1}$, we use the classifying space $BX$, which satisfies $H^*(BX;\F_p)=\F_p[y_4,y_{2p+2}]$. We will prove that
$\P^1(y_4)=y_{2p+2}+Ay_4^{(p+1)/2}$, for some generator $y_{2p+2}$ and some $A\in \F_p$, from which the desired result about the $x$'s follows
immediately from the map $\Sigma X\to BX$, which in $H^*(-;\F_p)$ sends $y_{j+1}$ to $x_j$ and sends products to 0.

First note that
\begin{eqnarray*}\P^1(y_4)&=&Ay_4^{(p+1)/2}+By_{2p+2}\\
\P^1(y_{2p+2})&=&Cy_4^p+Dy_4^{(p-1)/2}y_{2p+2},\end{eqnarray*} for some $A$, $B$, $C$, $D$ in $\F_p$. By the unstable property of the Steenrod
algebra,
\begin{equation}\label{uns}\P^{p+1}(y_{2p+2})=y_{2p+2}^p.\end{equation}
We must have
$$\P^p(y_{2p+2})=\sum_{j=0}^{p-1}c_jy_4^{1+j(p+1)/2}y_{2p+2}^{p-1-j},$$
for some $c_j\in \F_p$. Since $\P^{p+1}=\P^1\P^p$ and
\begin{equation}\label{Cartan}\P^1(y_4^iy_{2p+2}^j)=i\P^1(y_4)y_4^{i-1}y_{2p+2}^j+jy_4^i\P^1(y_{2p+2})y_{2p+2}^{j-1},
\end{equation}
the only way to obtain (\ref{uns}) is if $c_0B=1$ in $\F_p$. Thus $B$ must be a unit, and the generator $y_{2p+2}$ can be chosen so that
$B=1$.\end{pf}

\section{Nonmodular individual cases}\label{nonmodsec}
In this section, we consider all cases 4 through 34, excluding case 28 (which is $F_4$), in the Shephard-Todd numbering in which $p$ does not divide the order of the
reflection group. We obtain a very attractive result. One modular case, $(X_{12},p=3)$ is also included here. There is some overlap of our
methods and results here with those in \cite{Lev}.

\begin{thm} Let $X=(X_n)_p$ with $4\le n\le 34$ and $n\ne 28$,
 excluding
the modular cases $(X_{29})_5$, $(X_{31})_5$, and $(X_{34})_7$, which will be considered in the next two sections. Then $X\simeq\prod S^{2d-1}$, where
$2d$ ranges over the integers listed as the ``type" in \cite{CE},  except for the 31 cases listed in Table \ref{gtbl}. In these, each $B(-,\ldots,-)$ is
built by fibrations from spheres of the indicated dimensions, with $\a_1$ as each attaching map, and occurs as a factor in a product
decomposition of the $p$-completion of some $SU(N)$.\label{thm3}\end{thm}

We will call the integers $d$, which are 1/2 times the ``type " numbers of Clark-Ewing, the ``degrees."

\begin{center}
\begin{minipage}{6.5in}
\begin{tab} {\bf Cases in Theorem \ref{thm3} which are not products of spheres}
\begin{center}

\label{gtbl}
\begin{tabular}{c|c|l}
Case&Prime&Space\\
\hline
$5$&$7$&$B(11,23)$\\
$8$&$5$&$B(15,23)$\\
$9$&$17$&$B(15,47)$\\
$10$&$13$&$B(23,47)$\\
$12$&$3$&$B(11,15)$\\
$14$&$19$&$B(11,47)$\\
$16$&$11$&$B(39,59)$\\
$17$&$41$&$B(39,119)$\\
$18$&$31$&$B(59,119)$\\
$20$&$19$&$B(23,59)$\\
$24$&$11$&$B(7,27)\times S^{11}$\\
$25$&$7$&$B(11,23)\times S^{17}$\\
$26$&$7$&$B(11,23,35)$\\
$26$&$13$&$B(11,35)\times S^{23}$\\
$27$&$19$&$B(23,59)\times S^{13}$\\
$29$&$13$&$B(15,39)\times S^7\times S^{23}$\\
$29$&$17$&$B(7,39)\times S^{15}\times S^{23}$\\
$30$&$11$&$B(3,23)\times B(39,59)$\\
$30$&$19$&$B(3,39)\times S^{27}\times S^{59}$\\
$30$&$29$&$B(3,59)\times S^{23}\times S^{39}$\\
$31$&$13$&$B(15,39)\times B(23,47)$\\
$31$&$17$&$B(15,47)\times S^{23}\times S^{39}$\\
$32$&$7$&$B(23,35,47,59)$\\
$32$&$13$&$B(23,47)\times B(35,59)$\\
$32$&$19$&$B(23,59)\times S^{35}\times S^{47}$\\
$33$&$7$&$B(7,19)\times B(11,23,35)$\\
$33$&$13$&$B(11,35)\times S^7\times S^{19}\times S^{23}$\\
$34$&$13$&$B(11,35,59,83)\times B(23,47)$\\
$34$&$19$&$B(11,47,83)\times B(23,59)\times S^{35}$\\
$34$&$31$&$B(23,83)\times S^{11}\times S^{35}\times S^{47}\times S^{59}$\\
$34$&$37$&$B(11,83)\times S^{23}\times S^{35}\times S^{47}\times S^{59}$
\end{tabular}
\end{center}
\end{tab}
\end{minipage}
\end{center}

\begin{rmk}{\rm The $v_1$-periodic homotopy groups of $B(2n+1,2n+2p-1)$ were obtained
in \cite[1.3]{BDMi}. Those of $B(11,23,35)_7$ and $B(23,35,47,59)_7$ were obtained in \cite[1.4]{BDMi}. Using \cite[1.5,1.9]{DY}, we find that
for $\eps=0,1$,
$$\vp_{2t-\eps}(B(11,35,59,83))_{(13)}\approx \begin{cases}0&t\not\equiv5\quad(12)\\
\Z/13^{\max(f_5(t),f_{17}(t),f_{29}(t),f_{41}(t))}&t\equiv5\quad(12),\end{cases}$$ where $f_\g(t)=\min(\g,4+\nu_{13}(t-\g))$, while
$$\vp_{2t-\eps}(B(11,47,83))_{(19)}\approx\begin{cases}0&t\not\equiv5\quad (18)\\
\Z/19^{\max(f'_5(t),f'_{23}(t),f'_{41}(t))}&t\equiv5\quad (18),\end{cases}$$ where $f'_\g(t)=\min(\g,3+\nu_{19}(t-\g))$.} \label{rmk2}
\end{rmk}

\begin{pf*}{Proof of Theorem \ref{thm3}} It is straightforward to check that the pairs
(case, prime) listed in Table \ref{gtbl} are the only non-modular cases in \cite[Table 1]{And} in which an admissible prime $p$ satisfies
that $(p-1)$ divides the difference of distinct degrees. Indeed all other admissible primes have $(p-1)$ greater than the maximum difference of
degrees. For example, Case 30 requires $p\equiv 1,4$ mod 5, and the degrees are 2, 12, 20, 30. The first few  primes of the required congruence
are 11,  19, and 29. Clearly 10, 18, and 28 divide differences of these degrees, but no larger $(p-1)$ can. Thus the unstable Adams spectral
sequence argument used in proving Theorem \ref{m|p+1} works the same way here to show that $X$ is a product of $S^{2d-1}$ in all cases not
appearing in Table \ref{gtbl}.  In the relevant range, the $E_2$-term will consist only of infinite towers, one for each generator. The first
deviation from that is a $\Z/p$ in filtration 1 in homotopy dimension $(2d-1)+(2p-3)$, where $d$ is the smallest degree. This will always be
greater than the dimension of the largest $S^{2d-1}$.

The next step is to show that the Steenrod operation $\P^1$ in $H^*(X;\F_p)$ must connect all the classes listed as adjacent generators in
one of the $B$-spaces in Table \ref{gtbl}. This was achieved independently of, and slightly earlier than, the author in \cite{Nish} and \cite{HNO}.
We include several of our proofs, omitting the most complicated cases, to illustrate our methods and for the benefit of the reader without access
to \cite{Nish} and \cite{HNO}.
We accomplish this by
considering the $A$-module $H^*(BX;\F_p)$. With one exception\footnote{The exception is $(X_{34})_{19}$, which is handled in \cite{Nish}.
We thank Nishinobu for pointing out a gap in the argument for this case which appeared in an earlier version of this paper.},
all cases involving factors of $B(2m-1,2m+2p-3)$ are implied by Lemma \ref{B2} by applying
$H^*(BX)\to H^{*-1}(X)$, which sends products to 0. Similarly, Lemma \ref{3f} covers the two cases with a factor $B(11,23,35)$.

Now we must show that the spaces $X$ have the homotopy type claimed. The first 10 cases are immediate from \cite[7.1]{MNT}, and the two other
non-product cases, i.e. $(X_{26})_7$ and $(X_{32})_7$,  follow from \cite[7.2,7.6]{MNT}. Note that these results of \cite{MNT} did not deal with
$p$-completed spaces, but the obstruction theory arguments used there apply in the $p$-complete context.
 There are 15 additional types which we claim to be quasi $p$-regular.
As defined in \cite{MT}, a space is quasi $p$-regular if it is $p$-equivalent to a product of spheres and spaces of the form
 $B(2n+1,2n+2p-1)$. In \cite{MT}
(see esp. \cite[pp. 330-334]{MT}), many exceptional Lie groups are shown to be quasi $p$-regular (for appropriate $p$) using a skeletal
approach. We could use that approach here, but we prefer to use the unstable Adams spectral sequence (UASS). The two methods are really
equivalent.

Let $q=2p-2$. In Diagrams \ref{sphch} and \ref{Bch},  we illustrate the UASS for $S^{2n+1}$ in dimension less than $2n+pq-1$ and for
$B(2n+1,2n+q+1)$ in dimension less than $2n+3q-3$. Diagram \ref{sphch} gives a nice interpretation of the statement of the homotopy groups in
\cite[13.4]{Toda}. If $n\ge p$, the paired dots in Diagram \ref{sphch} will not occur in the pictured range. The nice thing about these charts
is that the $\F_p$-cohomology groups of our spaces $X$ are known to agree with that of their putative product decomposition as unstable algebras
over the Steenrod algebra, and are of the required universal form for the UASS to apply; hence their UASS has $E_2$-term the sum of the relevant
charts of spheres and $B$-spaces. In all cases, there will be no possible differentials.

One can check that in all 15 cases in which $X$ is claimed to be quasi $p$-regular, the towers in UASS$(X)$ corresponding to the spheres and the
bottom cell of each $B(2n+1,2n+q+1)$ cannot support a differential, and hence yield maps from the sphere or $S^{2n+1}$ into $X$. Next one checks
that $\pi_{2n+q}(X)=0$ and $\pi_{4n+q+1}(X)=0$. As these are the groups in which the obstruction to extending the map $S^{2n+1}\to X$ over
$B(2n+1,2n+q+1)$ lie, we obtain the desired extension. Finally, we take the product of maps $B\to X$ and $S^{2d_i-1}\to X$, using the group
structure of $X$, to obtain the desired $p$-equivalence from a product of spheres and $B$-spaces into $X$.

The remaining cases, $(X_{33})_7$, $(X_{34})_{13}$, and $(X_{34})_{19}$, are handled similarly. The $E_2$-term of the UASS converging to
$\pi_*(X)$ is isomorphic to that of its putative product decomposition. For example, $E_2(X_{34})_{13}$ is the sum of Diagram \ref{Bch} with
$n=11$ and $q=24$ plus Diagram \ref{3ch}. We can map $S^{23}\to X$ and $S^{11}\to X$ corresponding to generators of homotopy groups. Then we can
extend the first map over the 47- and 70-cells because $\pi_{46}(X)=0$ and $\pi_{69}(X)=0$. This gives a map $B(23,47)\to X$. Similarly we can
extend the second map over cells of $B(11,35,59,83)$ of dimension 46, 70, 94, 118, 142, 105, 129, 153, and 188. Taking the product of these two
maps, using the multiplication of $X$, yields the desired 13-equivalence $B(23,47)\times B(11,35,59,83)\to X$. The other two cases are handled
similarly.
\end{pf*}

\begin{center}
 \begin{minipage}{6.5in}
 \begin{diag}\label{sphch}{UASS$(S^{2n+1})$ in dim $<2n+pq-1$. Here $n<p$.}
\begin{center}
  \begin{picture}(480,210)(30,0)
  \def\mp{\multiput}
  \def\elt{\circle*{3}}
\def\lelt{\circle*{1}}
\put(0,20){\line(1,0){480}} \put(10,20){\vector(0,1){190}} \mp(10,20)(0,20){4}{\elt} \mp(70,40)(60,20){2}{\elt} \mp(190,85)(10,5){3}{\lelt}
\put(270,120){\elt} \mp(330,140)(10,0){2}{\elt} \mp(460,200)(10,0){2}{\elt} \mp(395,165)(10,5){3}{\lelt} \put(0,10){$\ssize{2n+1}$}
\put(60,10){$\ssize{2n+q}$} \put(120,10){$\ssize{2n+2q}$} \put(260,10){$\ssize{2n+nq}$} \put(320,10){$\ssize{2n+(n+1)q}$}
\put(450,10){$\ssize{2n+(p-1)q}$}
\end{picture}
\end{center}
\end{diag}
\end{minipage}
\end{center}

\noindent Where there is a pair of dots, the grading at the bottom refers to the one on right, and the other is in grading 1 less.

\begin{center}
\begin{minipage}{6.5in}
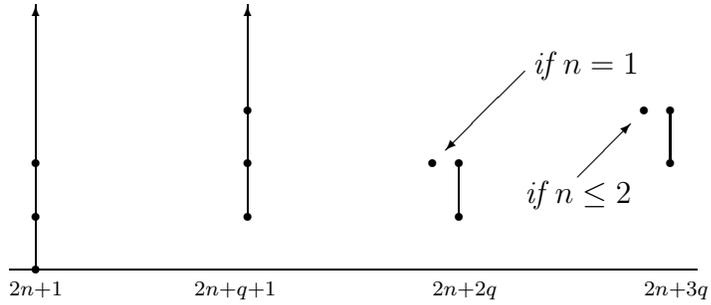
\begin{diag}\label{Bch}{UASS$(B(2n+1,2n+q+1))$ in dim $<2n+3q-3$.}
\begin{center}
\begin{picture}(260,130)
\put(0,20){\line(1,0){260}} \put(10,20){\vector(0,1){100}} \put(90,40){\vector(0,1){80}} \multiput(10,20)(0,20){3}{\circle*{3}}
\multiput(90,40)(0,20){3}{\circle*{3}} \put(170,40){\line(0,1){20}} \multiput(170,40)(0,20){2}{\circle*{3}} \put(160,60){\circle*{3}}
\put(195,95){\vector(-1,-1){30}} \put(250,60){\line(0,1){20}} \multiput(250,60)(0,20){2}{\circle*{3}} \put(240,80){\circle*{3}}
\put(215,55){\vector(1,1){20}} \put(195,45){$\text{if }n\le2$} \put(198,94){$\text{if }n=1$} \put(0,10){$\ssize{2n+1}$}
\put(70,10){$\ssize{2n+q+1}$} \put(160,10){$\ssize{2n+2q}$} \put(240,10){$\ssize{2n+3q}$}
\end{picture}
\end{center}
\end{diag}
\end{minipage}
\end{center}

\begin{center}
 \begin{minipage}{6.5in}
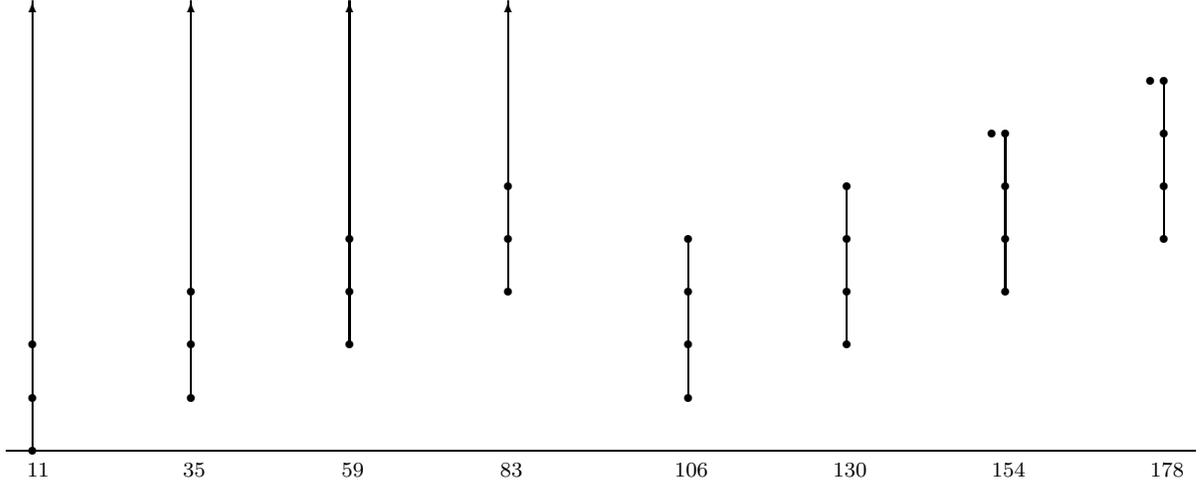
 \begin{diag}\label{3ch}{UASS$(B(11,35,59,83)_{13})$ in dim\ $<200$.}
\begin{center}
  \begin{picture}(480,200)(30,0)
  \def\mp{\multiput}
  \def\elt{\circle*{3}}
\put(0,20){\line(1,0){450}} \put(10,20){\vector(0,1){170}} \put(70,40){\vector(0,1){150}} \put(130,60){\vector(0,1){130}}
\put(190,80){\vector(0,1){110}} \mp(10,20)(0,20){3}{\elt} \mp(70,40)(0,20){3}{\elt} \mp(130,60)(0,20){3}{\elt} \mp(190,80)(0,20){3}{\elt}
\mp(258,40)(60,20){4}{\line(0,1){60}} \mp(258,40)(0,20){4}{\elt} \mp(318,60)(0,20){4}{\elt} \mp(378,80)(0,20){4}{\elt}
\mp(438,100)(0,20){4}{\elt} \mp(373,140)(60,20){2}{\elt} \put(8,10){$\ssize{11}$} \put(67,10){$\ssize{35}$} \put(127,10){$\ssize{59}$}
\put(187,10){$\ssize{83}$} \put(253,10){$\ssize{106}$} \put(313,10){$\ssize{130}$} \put(373,10){$\ssize{154}$} \put(433,10){$\ssize{178}$}
\end{picture}
\end{center}
\end{diag}
\end{minipage}
\end{center}

An alternate proof of Theorem \ref{thm3} can be obtained using \cite[1.3,1.4]{CHZ}, which can be interpreted as the following theorem.
In Section \ref{Harsec}, we will provide a proof of a strengthened version of this result.
\begin{thm}\label{CHZthm}$(\cite[1.3,1.4]{CHZ})$
If $X$ is an $H$-space of rank $r<p-1$ with torsion-free homology, then there are $H$-spaces $X_1,\ldots,X_r$ with $X_1=S^{n_1}$ and $X_r=X$,
and there are fibrations $X_{i-1}\to X_i\to S^{n_i}$ for $2\le i\le r$ as in the diagram

\begin{equation}\label{fibrs}\begin{CD} X_1@>>> X_2@>>> \cdots @>>> X_r\\
@. @VVV @. @VVV\\ @. S^{n_2} @. @.S^{n_r}.\end{CD}\end{equation} The homotopy type of the $p$-localization of $X$ is determined by the elements of
$\pi_{n_i-1}(X_{i-1})$ associated to these fibrations.\end{thm}

In order to apply this, one would still need to determine the $\P^1$-action and to check that the relevant homotopy groups
$\pi_{n_i-1}(X_{i-1})$ are cyclic.

In the following lemmas, which were used above, $g_i$ denotes a generator in grading $i$.
\begin{lem}\label{B2} a. If $m\not\equiv 1$ mod $p$ and $\F_p[g_{2m},g_{2m+2p-2}]$ is an unstable
$A$-algebra, then $\P^1g_{2m}\equiv ug_{2m+2p-2}$ mod decomposables, with $u\ne0$.

b. The same conclusion holds if the unstable $A$-algebra contains additional generators in dimensions $d\not\equiv 2m$ mod $(2p-2)$,
provided also that $d+2p-2\not\equiv0$ mod $(2m+2p-2)$.
\end{lem}
\begin{pf} a. For dimensional reasons, we must have $\P^1g_{2m}=\a g_{2m+2p-2}$ plus possibly a power
of $g_{2m}$, for some $\a\in \F_p$, and $\P^1g_{2m+2p-2}=g_{2m}Y$, for some polynomial $Y$. The unstable condition requires that
$\P^{m+p-1}g_{2m+2p-2}=g_{2m+2p-2}^p$, and, since $m\not\equiv1$ mod $p$, this equals, up to unit, $\P^1\P^{m+p-2}g_{2m+2p-2}$. For dimensional
reasons,
\begin{equation}\label{-1}\P^{m+p-2}g_{2m+2p-2}=\b g_{2m}g_{2m+2p-2}^{p-1}+g_{2m}^3Z\end{equation}
for some $\b\in \F_p$ and some polynomial $Z$. By the Cartan formula (similar to (\ref{Cartan})), the only way that $\P^1$ applied to (\ref{-1})
can yield $g_{2m+2p-2}^p$ is if both $\a$ and $\b$ are units.

b. One way that the additional generators could affect  the argument for part (a) would be if several of them (possibly the same
one) were multiplied together to get into the congruence of part (a). By the Cartan formula, $\P^1$ of such a product will still involve some of
these additional generators as factors, and so cannot yield the $g_{2m+2p-2}^p$ term on which the argument focuses.
The other way, pointed out to the author by Nishinobu, would be if there were a generator $g_d$ such that $\P^1g_d=u
g_{2m+2p-2}^i$ for some $u\ne0$ and $i>0$, an eventuality excluded by our second hypothesis. If there were such a $g_d$ and
also $\P^{m+p-2}g_{2m+2p-2}$ included a
term $g_dg_{2m+2p-2}^{p-i}$, this would provide an alternative way to achieve $g_{2m+2p-2}^p$ in $\P^{m+p-1}g_{2m+2p-2}$.  \end{pf}

\begin{lem}\label{3f}
 If $\F_7[g_{12},g_{24},g_{36}]$ is an unstable $A$-algebra, then, mod decomposables, $\P^1g_{12}=u_1g_{24}$ and $\P^1g_{24}=u_2g_{36}$ with
$u_i\ne0$. The same conclusion holds for $\F_7[g_{12},g_{24},g_{36},g_8,g_{20}]$.
\end{lem}
\begin{pf}
The first part is a direct consequence of the second, so we just consider $\F_7[g_{12},g_{24},g_{36},g_8,g_{20}]$.
  We work mod the ideal  $(g_{12},g_8,g_{20})$. Then $\P^1g_{12}\equiv\a g_{24}$, $\P^1g_{24}\equiv\b g_{36}$, and $\P^1g_{36}\equiv\g
g_{24}^2$, for some $\a$, $\b$, $\g$ in $\F_7$. This latter term complicates things somewhat. We also have $\P^1g_8\equiv0$ and
$\P^1g_{20}\equiv0$.

The unstable condition implies $\P^1\P^{17}g_{36}= ug_{36}^7$ with $u\ne0$. We use the Cartan formula as in the previous proof. The only way to get
to $g_{36}^7$ by $\P^1$  is if $\b\ne0$, implying the result for $\P^1g_{24}$.

However, there are two ways that $\P^1\P^{11}g_{24}$ might yield $g_{24}^7$, one via
$\P^1(g_{12}g_{24}^6)$ and the other via $\P^1(g_{24}^5g_{36})$. Instead, we consider $(\P^1)^5\P^7g_{24}$. We must have
$$\P^7g_{24}\equiv\delta_1g_{36}^3+\delta_2g_{24}^3g_{36}$$
for some $\delta_i\in \F_7$. We compute
\begin{eqnarray*}(\P^1)^5g_{36}^3&\equiv&\b\g^4g_{24}^7-\b^3\g^2g_{24}g_{36}^4\\
(\P^1)^5(g_{24}^3g_{36})&\equiv&\b^2\g^3g_{24}^7+5\b^3\g^2g_{24}^4g_{36}^2.\end{eqnarray*} Assuming that $\a=0$, so that the omitted terms of
the form $g_{12}Y$ in $\P^7g_{24}$ have $(\P^1)^5(\text{them})\equiv0$, then we obtain
\begin{eqnarray*}ug_{24}^7&\equiv& (\P^1)^5\P^7(g_{24})\\
&\equiv&(\P^1)^5(\delta_1g_{36}^3+\delta_2g_{24}^3g_{36})\\
&\equiv&\delta_1(\b\g^4g_{24}^7-\b^3\g^2g_{24}g_{36}^4)+\delta_2(\b^2\g^3g_{24}^7+5\b^3\g^2g_{24}^4g_{36}^2).
\end{eqnarray*}
Here we have also used that omitted terms divisible by $g_8$ or $g_{20}$ have $(\P^1)^5(\text{them})\equiv0$.
Coefficients of $g_{24}^7$ imply $\b\ne0$, $\g\ne0$, and some $\delta_i\ne0$, but this then gives a contradiction regarding $g_{24}g_{36}^4$ or
$g_{24}^4g_{36}^2$. Thus the assumption that $\a=0$ must have been false.
\end{pf}

\section{5-primary modular cases}\label{modsec}
In this section, we determine the $v_1$-periodic homotopy groups of the modular 5-compact groups $X_{29}$ and $X_{31}$. We pass directly from
invariant polynomials to Adams operations in $K^*(X)$ and thence to $\vp_*(X)$.  In
Theorems \ref{29conj} and \ref{X31E8}, we relate the homotopy type of $(X_{29})_5$
and $(X_{31})_5$ to that of $SU(20)$ and $E_8$. Theorem \ref{29conj} was conjectured by the author in an earlier version of this paper.
It and Theorem \ref{X31E8} and their proofs were provided to the author by John Harper.

The input to determining the Adams module $K^*(X_{29};\hat\Z_5)$ is the following result due to Aguad\'e (\cite{Ag}) and Maschke (\cite{Mas}).
Throughout the rest of the paper, we will denote by $m_{(e_1,\ldots,e_k)}$ the smallest symmetric polynomial on variables $x_1,\ldots,x_\ell$
(the value of $\ell$ will be implicit) containing the term $x_1^{e_1}\cdots x_k^{e_k}$.
\begin{thm} There is a reflection group $G_{29}$ acting on $(\hat\Z_5)^4$, and there is
a space $BX_{29}$ and map $BT\to BX_{29}$ with $BT=K((\hat\Z_5)^4,2)$ such that
$$H^*(BX_{29};\hat\Z_5)\approx H^*(BT;\hat\Z_5)^{G_{29}},$$
the invariants under the natural action of $G_{29}$ on $H^*(BT;\hat\Z_5)=\hat\Z_5[x_1,x_2,x_3,x_4]$ with $|x_i|=2$. Moreover,
$H^*(BT;\hat\Z_5)^{G_{29}}$ is a polynomial algebra on the following four invariant polynomials:
\begin{eqnarray*}f_4&=&m_{(4)}-12m_{(1,1,1,1)}\\
f_8&=&m_{(8)}+14m_{(4,4)}+168m_{(2,2,2,2)}\\
f_{12}&=&m_{(12)}-33m_{(8,4)}+330m_{(4,4,4)}+792m_{(6,2,2,2)}\\
f_{20}&=&m_{(20)}-19m_{(16,4)}-494m_{(12,8)}-336m_{(14,2,2,2)}+716m_{(12,4,4)}\\
&&+1038m_{(8,8,4)}+7632m_{(10,6,2,2)}+129012m_{(8,4,4,4)}+106848m_{(6,6,6,2)}.
\end{eqnarray*}\label{Masthm}
\end{thm}
\begin{pf} The group $G_{29}$ is the subgroup of $GL(\bC,4)$ generated by the following four matrices.
These can be seen explicitly in \cite{Ag}.
$$\frac12\begin{pmatrix}1&-1&-1&-1\\ -1&1&-1&-1\\ -1&-1&1&-1\\ -1&-1&-1&1\end{pmatrix},\qquad
\begin{pmatrix}0&i&0&0\\-i&0&0&0\\ 0&0&1&0\\ 0&0&0&1\end{pmatrix},\qquad \begin{pmatrix}0&1&0&0\\ 1&0&0&0\\
0&0&1&0\\ 0&0&0&1\end{pmatrix},\qquad \begin{pmatrix} 1&0&0&0\\ 0&0&1&0\\ 0&1&0&0\\ 0&0&0&1\end{pmatrix}.$$ Since $i\in\hat\Z_5$, these act on
$(\hat\Z_5)^4$, and this induces an action on $$H^*(K((\hat\Z_5)^4,2))\approx\hat\Z_5[x_1,x_2,x_3,x_4].$$ The invariants of this action were
determined by Maschke (\cite{Mas}) to be the polynomials stated in the theorem.  Although he did not state them all explicitly, they can be
easily generated by: (a) define $\phi$, $\psi_i$, and $\chi$ as on his page 501, then (b) define $\Phi_1,\ldots,\Phi_6$ as on his page 504, and
finally (c) let $f_4=-\frac12\Phi_6$ and $f_8=F_8$, $f_{12}=F_{12}$, and $f_{20}=F_{20}$ as on his page 505. See also \cite[p.287]{ST} for a
reference to this work.

Actually, Maschke's work and that of \cite{ST} involved finding generators for the complex invariant ring. To see that these
integral polynomials generate the invariant ring over $\hat\Z_5$, one must show that they cannot be decomposed over $\Z/5$.
For example, one must verify that $f_{20}$ cannot be decomposed mod 5 as a linear combination of $f_8f_{12}$, $f_4^2f_{12}$, $f_4f_8^2$,
$f_4^3f_8$, and $f_4^5$. The need to do this was pointed out to the author by Kasper Andersen in a dramatic way, as will be described
prior to \ref{Kgen}. The verification here was performed by Andersen using a {\tt Magma} program.

Aguad\'e (\cite{Ag}) constructed the 5-compact group $(BX,X)$ corresponding to this modular reflection group.
\end{pf}

The approach based on the following proposition benefits from a suggestion of Clarence Wilkerson.
\begin{prop} Let $(BX,X)$ be a $p$-compact group corresponding to a reflection group $G$ acting on $BT=K((\hat\Z_p)^n,2)$. Suppose
$H^*(BT;\hat\bQ_p)^G=\hat\bQ_p[f_1,\ldots,f_k]$, where $f_i$ is a polynomial in $y_1,\ldots,y_n$ with $y_j\in H^2(BT;\hat\bQ_p)$
corresponding to the $j^{\text{th}}$ factor. Let $K^*(BT;\hat\bQ_p)=\hat\bQ_p[\![x_1,\ldots,x_n]\!]$ with $x_i$ the class of $H-1$ in the
$i^{\text{th}}$ factor,
where $H$ is the complex Hopf bundle. Let $\ell_0(x)=\log(1+x)$. Then
$$K^*(BX;\hat\Z_p)\approx \hat\bQ_p[\![f_1(\ell_0(x_1),\ldots,\ell_0(x_n)),\ldots,f_k(\ell_0(x_1),\ldots,\ell_0(x_n))]\!]\cap
\hat\Z_p[\![x_1,\ldots,x_n]\!].$$
\end{prop}
\begin{pf} The Chern character $K^*(BT;\hat\bQ_p)@>\text{ch}>> H^*(BT;\hat\bQ_p)$ satisfies $\text{ch}(\ell_0(x_i))=y_i$ and hence,
since ch is a ring homomorphism, $\text{ch}(f_j(\el_0(x_1),\ldots,\el_0(x_n)))=f_j(y_1,\ldots,y_n)$. It commutes with the action of $G$, and
hence sends invariants to invariants. Indeed \begin{equation}\label{Qp}
K^*(BT;\hat\bQ_p)^G=\hat\bQ_p[\![f_1(\ell_0(x_1),\ldots,\ell_0(x_n)),\ldots,f_k(\ell_0(x_1),\ldots,\ell_0(x_n))]\!].
\end{equation}
The invariant ring in $K^*(BT;\hat\Z_p)$ is just the intersection of (\ref{Qp}) with $\hat\Z_p[\![x_1,\ldots,x_n]\!]$. Finally we use a result
of \cite{JO} that $K^*(BX;\hat\Z_p)\approx K^*(BT;\hat\Z_p)^G$.\end{pf}

Thus with $f_4,f_8,f_{12},f_{20}$  as in \ref{Masthm}, we wish to find algebraic combinations of $$f_4(\el_0(x_1),\ldots,\el_0(x_4)),
\ldots,f_{20}(\el_0(x_1),\ldots,\el_0(x_4))$$
 which have coefficients in $\hat\Z_5$. A theorem of \cite{JO} which states that for a
$p$-compact group $BX$ there is an isomorphism $K^*(BX;\hat\Z_p)\approx\hat\Z_p[\![g_1,\ldots,g_k]\!]$, and the collapsing, for dimensional
reasons, of the Atiyah-Hirzebruch spectral sequence
\begin{equation}\label{AHSS}H^*(BX;K^*(\text{pt};\hat\Z_p))\Rightarrow K^*(BX;\hat\Z_p)\end{equation}
implies that the generators $g_j$ can be chosen to be of the form $f_j(x_1,\ldots,x_n)$ mod higher degree polynomials.

Finding these algebraic combinations can be facilitated by using the $p$-typical log series
$$\el_p(x)=\sum_{i\ge0}x^{p^n}/p^n.$$
By \cite{Has}, there is a series $h(x)\in\Z_{(p)}[\![x]\!]$ such that $\el_0(h(x))=\el_p(x)$ and $h(x)\equiv x$ mod $(x^2)$. Let $x_i'=h(x_i)$.
For any $c_{\bold e}\in \hat\bQ_p$ with ${\bold e}=(e_4,e_8,e_{12},e_{20})$, we have
\begin{eqnarray}&&\sum c_{\bold e} f_4(\el_p(x_1),\ldots,\el_p(x_4))^{e_4}\cdots f_{20}(\el_p(x_1),\ldots,\el_p(x_4))^{e_{20}}\label{LHS}\\
&=&\sum c_{\bold e} f_4(\el_0(x'_1),\ldots, \el_0(x'_4))^{e_4}\cdots f_{20}(\el_0(x'_1),\ldots,\el_0(x'_4))^{e_{20}},\label{equal}\end{eqnarray}
where the sums are taken over various $\bold e$.
We will find $c_{\bold e}$ so that (\ref{LHS}) is in $\hat\Z_p[\![x_1,\ldots,x_4]\!]$. Thus so is (\ref{equal}), and hence also $\sum
c_{\bold e} f_4(\el_0(x_1),\ldots,\el_0(x_4))^{e_4}\cdots f_{20}(\el_0(x_1),\ldots,\el_0(x_4))^{e_{20}}$, since $h(x)\in\Z_{(p)}[\![x]\!]$.

A {\tt Maple} program, which will be described in the proof, was used to prove the following result.
\begin{thm} \label{29Kpolys} Let $f_4,f_8,f_{12},f_{20}$ be as in \ref{Masthm}, and let
$$F_j=F_j(x_1,\ldots,x_4)=f_j(\el_0(x_1),\ldots,\el_0(x_4)).$$
Then the following  series are $5$-integral through grading 20; i.e., their coefficients of all monomials $x_1^{e_1}\cdots x_4^{e_4}$ with $\sum e_i\le20$
are 5-integral.
\begin{eqnarray*}&&F_4-\tfrac1{10}F_4^2-\tfrac15F_8-\tfrac{16}{25}F_{12}-\tfrac7{25}F_4F_8+\tfrac4{25}F_4^3-\tfrac{13}{125}F_4F_{12}
-\tfrac{57}{125}F_4^2F_8\\ &&\qquad-\tfrac{102}{125}F_4^4-\tfrac{62}{125}F_8^2-\tfrac{64}{125}F_{20}-\tfrac4{625}F_4^5
-\tfrac{42}{125}F_4^2F_{12}-\tfrac{11}{25}F_4F_8^2-\tfrac{72}{125}F_8F_{12};\\
&&F_8-\tfrac85F_{12}-\tfrac7{25}F_8^2-\tfrac4{25}F_{20}-\tfrac{21}{125}F_8F_{12};\\
&&F_{12}-\tfrac25F_8^2-\tfrac15F_{20}-\tfrac4{25}F_8F_{12}.\end{eqnarray*} \end{thm}
\begin{pf} As observed in the paragraph preceding the theorem, it suffices to show that the same is true for
$\widetilde F_j=f_j(x_1+\frac15x_1^5,\ldots,x_4+\frac15x_4^5)$. The advantage of this is to decrease the number of terms which must be kept
track of and looked at. We work one grading at a time, expanding relevant products of $F_j$'s as combinations of monomial symmetric polynomials in the
fixed grading, and then solving a system of linear equations to find the combinations that work. We illustrate with the calculation for
modifications of $F_4$ in gradings 8 and then 12.

In grading 8, we have
\begin{eqnarray*} \Ft_4&=&\tfrac45m_{(8)}-\tfrac{12}5m_{(5,1,1,1)}\\
\Ft_8&=&m_{(8)}\qquad\qquad\quad+14m_{(4,4)}+168m_{(2,2,2,2)}\\
\Ft_4^2&=&m_{(8)}-24m_{(5,1,1,1)}+2m_{(4,4)}+144m_{(2,2,2,2)}.\end{eqnarray*} We wish to choose $a$ and $b$ so that, in grading 8, $\tfrac a5\Ft_8+\frac b5
\Ft_4^2\equiv \Ft_4$ mod integers. Thus we must solve a system of mod 5 equations for $a$ and $b$ with augmented matrix
$$\left(\begin{array}{cc|c} 1&1&4\\ 0&-24&-12\\ 14&2&0\\ 168&144&0\end{array}\right).$$
The solution is $a=1$, $b=1/2$. We could also have used $b=3$ since we are working mod 5.

Let $\Ft'_4=\Ft_4-\frac1{10}\Ft_4^2-\frac15\Ft_8$. In grading 12, we have
\begin{eqnarray*} \Ft'_4&=&-\tfrac6{25}m_{(12)}-\tfrac{12}5m_{(8,4)}\qquad\quad-\tfrac{96}5m_{(6,2,2,2)}+\tfrac{12}5m_{(9,1,1,1)}+\tfrac{12}{25}m_{(5,5,1,1)}\\
\Ft_{12}&=&m_{(12)}-33m_{(8,4)}+330m_{(4,4,4)}+792m_{(6,2,2,2)}\\
\Ft_8\Ft_{4}&=&m_{12}+15m_{(8,4)}+42m_{(4,4,4)}+168m_{(6,2,2,2)}-12m_{(9,1,1,1)}-168m_{(5,5,1,1)}-2016m_{(3,3,3,3)}\\
\Ft_4^3&=&m_{(12)}+3m_{(8,4)}+6m_{(4,4,4)}+432m_{(6,2,2,2)}-36m_{(9,1,1,1)}-72m_{(5,5,1,1)}-1728m_{(3,3,3,3)}.\end{eqnarray*} We wish to choose
$a$, $b$, and $c$ so that, in grading 12, $\tfrac a{25}\Ft_{12}+\frac b{25}\Ft_8\Ft_4+\frac c{25}\Ft_4^3\equiv \Ft'_4$ mod integers. Thus we must solve a system
of equations mod 25 whose augmented matrix is
$$\left(\begin{array}{ccc|c}1&1&1&-6\\ -33&15&3&-60\\ 330&42&6&0\\ 792&168&432&-480\\ 0&-12&-36&60\\ 0&-168&-72&12\\ 0&-2016&-1728&0\end{array}\right).$$
The solution is $a=16$, $b=7$, and $c=-4$.

We perform similar calculations for $\Ft_8$ in grading 12, then for $\Ft''_4$, $\Ft'_8$, and $\Ft_{12}$ in gradings 16 and then 20.
\end{pf}

By the observation in the paragraph involving (\ref{AHSS}), the modified versions of $F_4$, $F_8$, and $F_{12}$ given in Theorem \ref{29Kpolys},
and also $F_{20}$, can be modified similarly in all subsequent gradings, yielding generators of the power series algebra $K^*(BX_{29};\hat\Z_5)$
which we will call $G_4$, $G_8$, $G_{12}$, and $G_{20}$. By \cite{JO}, $K^*(X_{29};\hat\Z_5)$ is an exterior algebra on classes $z_3$, $z_7$,
$z_{11}$, and $z_{19}$ in $K^1(-)$ obtained using the map $e:\Sigma X=\Sigma\Omega BX\to BX$ and Bott periodicity $B:K^1(X)\to K^{-1}(X)$ by
$z_i=B^{-1}e^*(G_{i+1})$. The following determination of the Adams operations is essential for our work on $v_1$-periodic homotopy groups.
Here and elsewhere $QK^1(-)$ denotes the indecomposable quotient.

\begin{thm}\label{psik} The Adams operation $\psi^k$ in $QK^1(X_{29};\hat\Z_5)$ on the generators $z_3$, $z_7$,
$z_{11}$, and $z_{19}$ is given by the matrix
$$\begin{pmatrix}k^3&0&0&0\\[3pt] \frac15k^3-\frac15k^7&k^7&0&0\\[3pt]
\frac{24}{25}k^3-\frac8{25}k^7-\frac{16}{25}k^{11}&\frac85k^7-\frac85k^{11}&k^{11}&0\\[3pt]
\frac{92}{125}k^3-\frac{12}{125}k^7-\frac{16}{125}k^{11}-\frac{64}{125}k^{19}&\frac{12}{25}k^7-\frac8{25}k^{11}-\frac4{25}k^{19}&\frac15k^{11}-\frac15k^{19}&k^{19}
\end{pmatrix}.$$
\end{thm}
\begin{pf} We first note that $$\psi^k(\el_0(x))=\el_0(\psi^kx)=\el_0((x+1)^k-1)=\log((x+1)^k)=k\log(x+1)=k\el_0(x).$$
Since $F_{4j}$ is homogeneous of degree $4j$ in $\el_0(x_i)$, $\psi^k(F_{4j})=k^{4j}F_{4j}$. We can use this to determine $\psi^k$ on the generators
$G_i$ which are defined as algebraic combinations of $F_{4j}$'s. We then apply $e^*$ to this formula to obtain $\psi^k$ in
$K^{-1}(X_{29};\hat\Z_5)$. Since $e^*$ annihilates decomposables, we need consider only the linear terms in the expressions which express $G_i$
in terms of $F_{4j}$'s. On the basis (over $\hat\bQ_5$)  $\langle e^*(F_4),e^*(F_8),e^*(F_{12}),e^*(F_{20})\rangle$, the matrix of $\psi^k$ is
$D=\text{diag}(k^4,k^8,k^{12},k^{20})$. On the basis (over $\hat\Z_5$)  $$\langle e^*(G_4),e^*(G_8),e^*(G_{12}),e^*(G_{20})\rangle,$$ it is
$P^{-1}DP$, where
$$P=\begin{pmatrix}1&0&0&0\\[3pt] -\frac15&1&0&0\\[3pt] -\frac{16}{25}&-\frac85&1&0\\[3pt] -\frac{64}{125}&-\frac4{25}&-\frac15&1\end{pmatrix}$$
is the change-of-basis matrix, obtained using the linear terms in \ref{29Kpolys}.
The matrix in the statement of the theorem is obtained by dividing $P^{-1}DP$ by $k$, since $\psi^k$ in $K^1(-)$ corresponds to $\psi^k/k$ in $K^{-1}(-)$.
\end{pf}

We can use Theorem \ref{psik} to obtain the $v_1$-periodic homotopy groups of $(X_{29})_5$ as follows.
\begin{thm}\label{v*G29} The groups $\vp_*(X_{29})_{(5)}$ are given by
$$\vp_{2t-1}(X_{29})\approx\vp_{2t}(X_{29})\approx\begin{cases}0&t\not\equiv3\quad (4)\\
\Z/5^3&t\equiv3,15\quad (20)\\
\Z/5^{\min(8,3+\nu_5(t-7-4\cdot5^4))}&t\equiv7\quad(20)\\
\Z/5^{\min(12,3+\nu_5(t-11-4\cdot5^8))}&t\equiv11\quad (20)\\
\Z/5^{\min(20,3+\nu_5(t-19-12\cdot5^{16}))}&t\equiv19\quad(20).\end{cases}$$\end{thm}
\begin{pf} We use the result of \cite{Bo} that $\vp_{2t}(X)_{(5)}$ is presented by the matrix $\begin{pmatrix}(\Psi^5)^T\\ (\Psi^2)^T-2^tI\end{pmatrix}$.
We form this matrix by letting $k=5$ and 2 in the matrix of Theorem \ref{psik} and letting $x=2^t$, obtaining
\begin{equation}\label{psimat}\begin{pmatrix}125&-15600&-31274880&-9765631257408\\ 0&78125&-78000000&-3051773400000\\ 0&0&48828125&-3814687500000\\
0&0&0&19073486328125\\ 8-x&-24&-1344&-268704\\
0&128-x&-3072&-84480\\
0&0&2048-x&-104448\\
0&0&0&524288-x\end{pmatrix}.\end{equation} Pivoting on the units (over $\Z_{(5)}$) in positions (5,2) and (7,4) and removing their rows and
columns does not change the group presented. We now have a 6-by-2 matrix, whose nonzero entries are polynomials in $x$ of degree 1 or 2. If
$x\not\equiv3$ mod 5, which is equivalent to $t\not\equiv3$ mod 4, the bottom two rows are $\left(\begin{smallmatrix} u_1&B\\
0&u_2\end{smallmatrix}\right)$ with $u_i$ units, and so the group presented is 0.

Henceforth, we assume $x\equiv 3$ mod 5. The polynomial in new position (5,2) is nonzero mod 5 for such $x$, and so we pivot on it, and remove its
row and column. The five remaining entries are ratios of polynomials with denominator nonzero mod 5. Let $p_1,\ldots,p_5$ denote the polynomials
in the numerators. The group $\vp_{2t}(X_{29})_{(5)}$ is $\Z/5^e$, where $e=\min(\nu(p_1(x)),\ldots,\nu(p_5(x)))$, where $x=2^t$. We abbreviate
$\nu_5(-)$ to $\nu(-)$ throughout the remainder of this section. We have
\begin{eqnarray*}p_1&=&-71122941747658752+9480741773824512x-74067383851199x^2+33908441866x^3\\
p_2&=&-66750692556800000+8897903174800000x-69512640100000x^2+31789306250x^3\\
p_3&=&-8327872\cdot10^{10}+11101145\cdot10^9x-86731015625000x^2+39736328125x^3\\
p_4&=&4\cdot10^{19}-533203125\cdot10^{10}x+41656494140625000x^2-19073486328125x^3\\
p_5&=&1099511627776-146567856128x+1145324544x^2-526472x^3+x^4.\end{eqnarray*}

For values of $m$ listed in the table, we compute and present in Table \ref{exptbl} the tuples $(e_0,e_1,e_2,e_3)$ so that, up to units,
\begin{equation}\label{peq}p_i(2^m+y)=5^{e_0}+5^{e_1}y+5^{e_2}y^2+5^{e_3}y^3\end{equation} (plus $y^4$ if $i=5$). Considerable preliminary
calculation underlies the choice of these values of $m$.
\begin{tab}{Exponents of polynomials}
\label{exptbl}
\begin{center}
\begin{tabular}{c|ccccc}
&&&$i$&&\\
$m$&$1$&$2$&$3$&$4$&$5$\\
\hline $3$&$3,2,1,0$&$\infty,7,6,5$&$\infty,12,12,10$&$\infty,21,20,19$&$\infty,3,2,1$\\
 $15$&$3,2,1,0$&$8,7,6,5$&$13,12,11,10$&$22,21,20,19$&$4,4,3,2$\\
 $7+4\cdot5^4$&$8,2,1,0$&$8,7,6,5$&$17,12,11,10$&$26,21,21,19$&$8,3,2,1$\\
 $11+4\cdot5^8$&$12,2,1,0$&$12,7,7,5$&$13,12,11,10$&$\infty,21,20,19$&$12,3,2,1$\\
 $19+12\cdot5^{16}$&$23,2,2,0$&$20,7,6,5$&$21,12,11,10$&$22,21,20,19$&$20,3,2,1$
\end{tabular}
\end{center}
\end{tab}

Recall that $\nu(2^{4\cdot5^i}-1)=i+1$, as is easily proved by induction. Thus
\begin{equation}\label{2025}p(2^{m+20j})=p(2^m+2^m(2^{20j}-1))=p(2^m+25j\cdot u),\end{equation} with $u$ a unit. Hence
$$\min\{\nu(p_i(2^{3+20j})):1\le i\le5\}=3$$
since $$p_1(2^{3+20j})=p_1(2^3+5^2ju)=5^3+5^2\cdot5^2ju+5(5^2ju)^2+(5^2ju)^3,$$ omitting some unit coefficients. Here we have set $y=5^2ju$ in
(\ref{peq}). Replacing $3$ by $15$ yields an identical argument. This yields the second line of Theorem \ref{v*G29}.

We use Table \ref{exptbl} to show \begin{equation}\label{case5}\min\{\nu(p_i(2^{19+12\cdot5^{16}+20j})):1\le
i\le5\}=\min(20,4+\nu(j))=\min(20,3+\nu(20j)).\end{equation} Indeed, for $\nu(j)\le16$, the minimum is achieved when $i=1$, with the 4 coming as $2+2$ with one 2
being from the 25 in (\ref{2025}) and the other 2 being the first 2 in the last row of Table \ref{exptbl}. If $\nu(j)>16$, the minimum is
achieved when $i=2$, using the first 20 in the last row of \ref{exptbl}. The last case of Theorem \ref{v*G29} follows easily from (\ref{case5}),
and the other two parts of \ref{v*G29} are obtained similarly.

To see that $\vp_{2t-1}(X_{29})\approx\vp_{2t}(X_{29})$, we argue in three steps. First, the two groups have the same order using \cite[8.5]{Bo}
and the fact that the kernel and cokernel of an endomorphism of a finite group have equal orders. Second, by \cite[4.4]{DHHA}, a presentation of
$\vp_{2t-1}(X_{29})$ is given by $\begin{pmatrix} \Psi^5\\ \Psi^2-2^t\end{pmatrix}$, i.e. like that for $\vp_{2t}(X_{29})$ except that the two
submatrices are not transposed. Third, we pivot on this matrix, which is (\ref{psimat}) with the top and bottom transposed, and find that we
can pivot on units three times, so that the group presented is cyclic.
\end{pf}

One of the factors in the  product decomposition of $SU(20)_5$ given in \cite{MNT} is an $H$-space $B_3^5(5)$ whose $\F_5$-cohomology is an
exterior algebra on classes of grading 7, 15, 23, 31, and 39, and which is built from spheres of these dimensions by fibrations.  By \cite{Ya},
there is a product decomposition
$$(SU(20)/SU(15))_5\simeq S^{31}\times S^{33}\times S^{35}\times S^{37}\times S^{39}.$$
Let $B(7,15,23,39)$ denote the fiber of the composite
$$B_3^5(5)\to SU(20)_5\to (SU(20)/SU(15))_5@>\rho >>(S^{31})_5.$$
\begin{thm}\label{29conj} (Harper) There is a homotopy equivalence
$$(X_{29})_5\simeq B(7,15,23,39).$$
\end{thm}

Note that this result requires more than \ref{CHZthm} because the ranks of these $H$-spaces are not less than $p-1$.
We will provide Harper's proof of this result in Section \ref{Harsec}. Here we just remark that our work above is required in the proof,
for the entry in position $(4,3)$ of the matrix of \ref{psik} implies that the 39-cell of $(X_{29})_5$ is attached to the 23-cell by
$\a_2$, which is not detected by primary Steenrod operations. This information is required in order to compare the two spaces in \ref{29conj}.

We can determine the Adams operations and $v_1$-periodic homotopy groups of $(X_{31})_5$ by an argument very similar to that used above for
$(X_{29})_5$. We shall merely sketch. The analogue of Theorem \ref{Masthm} is
\begin{thm} There is an isomorphism $H^*(BX_{31};\hat\Z_5)\approx H^*(BT;\hat\Z_5)^{G_{31}}$, where $G_{31}$ has the four generators
given for $G_{29}$ in the proof of \ref{Masthm} and also $\left(\begin{smallmatrix}1&0&0&0\\ 0&1&0&0\\
0&0&-1&0\\0&0&0&1\end{smallmatrix}\right)$. Then $H^*(BT;\hat\Z_5)^{G_{31}}$ is a polynomial ring on the generators $f_8$, $f_{12}$, and
$f_{20}$ given in \ref{Masthm} together with
\begin{eqnarray*} f_{24}&=&m_{(24)}-66m_{(20,4)}+1023m_{(16,8)}+2180m_{(12,12)}+1293156m_{(8,8,4,4)}\\
&&+267096m_{(12,4,4,4)}
+2121984m_{(6,6,6,6)}+620352m_{(10,6,6,2)}-4032m_{(14,6,2,2)}\\
&&-190080m_{(10,10,2,2)}
-11892m_{(12,8,4)}-4938m_{(16,4,4)}-24534m_{(8,8,8)}\\
&&-2304m_{(18,2,2,2)}.\end{eqnarray*}\label{31polys}
\end{thm}

The analogue of Theorem \ref{29Kpolys} is
\begin{thm} Let $f_8,f_{12},f_{20},f_{24}$ be as in \ref{31polys}, and let
$$F_j=f_j(\ell_0(x_1),\ldots,\ell_0(x_4)).$$
Then the following series are 5-integral through grading 24.
\begin{eqnarray*}&&F_8-\tfrac85F_{12}-\tfrac7{25}F_8^2-\tfrac4{25}F_{20}-\tfrac{21}{125}F_8F_{12}-\tfrac{99}{125}F_{24}-\tfrac{597}{625}F_8^3-\tfrac{558}{625}F_{12}^2\\
&&F_{12}-\tfrac25F_8^2-\tfrac15F_{20}-\tfrac4{25}F_8F_{12}-\tfrac{18}{25}F_{24}-\tfrac{74}{125}F_8^3-\tfrac{11}{125}F_{12}^2\\
&&F_{20}-\tfrac35F_{24}-\tfrac25F_8^3.\end{eqnarray*}
\end{thm}

The analogue of \ref{psik} is
\begin{thm}\label{31ops} The Adams operation $\psi^k$ in $K^1(X_{31};\hat\Z_5)$ on the generators $z_7$, $z_{11}$, $z_{19}$, and
$z_{23}$ is given by the matrix
$$\begin{pmatrix}k^7&0&0&0\\[4pt]
 \tfrac85k^7-\tfrac85k^{11}&k^{11}&0&0\\[4pt]
  \tfrac{12}{25}k^7-\tfrac8{25}k^{11}-\tfrac4{25}k^{19}&\tfrac15k^{11}-\tfrac15k^{19}&k^{19}&0\\[4pt]
\tfrac{279}{125}k^7-\tfrac{168}{125}k^{11}-\tfrac{12}{125}k^{19}-\tfrac{99}{125}k^{23}&\tfrac{21}{25}k^{11}
-\tfrac3{25}k^{19}-\tfrac{18}{25}k^{23}&\tfrac35k^{19}-\tfrac35k^{23}&k^{23}\end{pmatrix}.$$
\end{thm}

The analogue of Theorem \ref{v*G29} is
\begin{thm} The groups $\vp_*(X_{31})_{(5)}$ are given by
$$\vp_{2t-1}(X_{31})\approx\vp_{2t}(X_{31})\approx\begin{cases}0&t\not\equiv3\quad (4)\\
\Z/5^3&t\equiv7,15\quad (20)\\
\Z/5^{\min(8,3+\nu_5(t-11-8\cdot5^4))}&t\equiv11\quad(20)\\
\Z/5^{\min(12,3+\nu_5(t-19-16\cdot5^8))}&t\equiv19\quad (20)\\
\Z/5^{\min(20,3+\nu_5(t-23-16\cdot5^{16}))}&t\equiv23\quad(20).\end{cases}$$\end{thm}

In \cite[Proposition 2.3]{Wilk}, it is proved that the exceptional Lie group $E_8$, localized at 5, admits a product decomposition as
$X_0(E_8)\times X_2(E_8)$ with $H^*(X_0(E_8);\F_5)\approx\Lambda(x_{15},x_{23},x_{39},x_{47})$. In Section \ref{Harsec},
we will prove the following result, which was pointed out by John Harper.
\begin{thm}\label{X31E8} (Harper) There is a homotopy equivalence
$$(X_{31})_5\simeq X_0(E_8).$$
\end{thm}

\section{The 7-primary modular case}\label{X34sec}
In this section, we first give in Theorem \ref{polys} new explicit formulas for the six polynomials which generate
as a polynomial algebra the invariant ring of the complex reflection group $G_{34}$ of \cite{ST}, called the Mitchell
group in \cite{CS}. Over $\Zh$, the invariant ring of $G_{34}$ is also a polynomial algebra, but the generators must
be altered slightly from the complex case, as we show prior to \ref{Kgen}. Next we use this information to find explicit generators
for $K^*(BX_{34};\Zh)$ in \ref{Kgen}, and from this the Adams operations in $QK^1(X_{34};\Zh)$ in \ref{Adams}.
These in turn enable us to compute the $v_1$-periodic homotopy groups $\vp_*(X_{34})_{(7)}$.
Finally, we prove in \ref{34conj} that $(X_{34})_7$ has the homotopy type of a space formed from $SU(42)$.
This result was conjectured by the author and proved by John Harper.

\begin{thm}\label{polys} The complex invariants of the reflection group $G_{34}$ (defined in the proof) form a polynomial algebra
$$\C[x_1,\ldots,x_6]^{G_{34}}\approx \C[f_6,f_{12},f_{18},f_{24},f_{30},f_{42}]$$ with generators given by
$$f_{6k}=(1+(-1)^k27^{k-1}\cdot5)m_{(6k)}+\sum_{s=1}^k\tbinom{6k}{3s}(1+(-1)^{k+s}27^{k-1})m_{(6k-3s,3s)}
+\sum_{\ebar}(\ebar)m_{\ebar},$$
where $\ebar$ ranges over all partitions $\ebar=(e_1,\ldots,e_r)$ of $6k$ with $3\le r\le6$ satisfying $e_i\equiv e_j$ mod 3
for all $i,j$, and $e_i\equiv 0$ mod 3 if $r<6$. Here also $(\ebar)$ denotes the multinomial coefficient
$(e_1+\cdots+e_r)!/(e_1!\cdots e_r!)$, and $m_\ebar$ the monomial symmetric polynomial, which is the shortest
symmetric polynomial in $x_1,\ldots,x_6$ containing $x_1^{e_1}\cdots x_r^{e_r}$.
\end{thm}
 For example, we have
\begin{itemize}
\item $f_6=-4m_{(6)}+40m_{3,3}+720m_{(1,1,1,1,1,1)}$;
\item $f_{12}=136m_{(12)}-26\binom{12}3m_{(9,3)}+28\binom{12}6m_{(6,6)}+\sum(\ebar)m_\ebar$, where $\ebar$ ranges
over
$$\{(6,3,3),\ (3,3,3,3),\ (2,2,2,2,2,2),\ (7,1,1,1,1,1),\ (4,4,1,1,1,1)\}.$$
\item $f_{18}=(1-5\cdot27^2)m_{(18)}+\binom{18}3(1+27^2)m_{(15,3)}+\binom{18}6(1-27^2)m_{(12,6)}+\binom{18}9(1+27^2)m_{(9,9)}
+\sum(\ebar)m_\ebar$, where $\ebar$ ranges
over
\begin{eqnarray*}&&\{(12,3,3),\ (9,6,3),\ (9,3,3,3),\ (6,6,6),\ (6,6,3,3),\ (6,3,3,3,3),\ (3,3,3,3,3,3),\\ &&(13,1,1,1,1,1),\ (10,4,1,1,1,1),\
(7,4,4,1,1,1),\ (4,4,4,4,1,1),\ (7,7,1,1,1,1),\\ &&(8,2,2,2,2,2),\ (5,5,2,2,2,2)\}
\end{eqnarray*}
\end{itemize}

\begin{pf*}{Proof of Theorem \ref{polys}}
As described in \cite{ST}, the reflection group $G_{34}$ is generated by reflections across the following hyperplanes in $\C^6$: $x_i-x_j=0$, $x_1-\omega
x_2=0$, and $x_1+x_2+x_3+x_4+x_5+x_6=0$. Here $\omega=e^{2\pi i/3}$. It follows easily that $G_{34}$ is generated by all permutation matrices
together with the following two:
\begin{equation}\label{twomatrices}\begin{pmatrix}0&\omega^2&0&0&0&0\\ \omega&0&0&0&0&0\\ 0&0&1&0&0&0\\ 0&0&0&1&0&0\\ 0&0&0&0&1&0\\
0&0&0&0&0&1\end{pmatrix},\qquad I-\tfrac13\begin{pmatrix}1&1&1&1&1&1\\1&1&1&1&1&1\\1&1&1&1&1&1\\1&1&1&1&1&1\\
1&1&1&1&1&1\\1&1&1&1&1&1\end{pmatrix}\end{equation}

In \cite{CS}, Conway and Sloane consider $G_{34}$ instead as the automorphisms of a certain $\Z[\omega]$-lattice in $\C^6$. The lattice has 756 vectors of
norm 2. There are none of smaller positive norm. 270 of these vectors are those with $\om^a$ in one position, $-\om^b$ in another, and 0 in the
rest. Here, of course, $a$ and $b$ can be 0, 1, or 2. The other 486 are those of the form $\pm\frac1{\sqrt{-3}}(\om^{a_1},\ldots,\om^{a_6})$
such that $\sum a_i\equiv 0$ mod 3.

As a partial verification that this lattice approach to $G_{34}$ is consistent with the reflection approach, one can verify that the reflection
matrices permute these 756 vectors.  It is obvious that permutation matrices do, and easily verified for the first matrix of
(\ref{twomatrices}). The second matrix of (\ref{twomatrices}), which has order 2, sends
\begin{itemize} \item$(\om,\om^2,0,0,0,0)$ to $\frac1{\sqrt{-3}}(\om^2,\om,1,1,1,1)$;
\item $\frac1{\sqrt{-3}}(1,1,1,1,1,1)$ to $-\frac1{\sqrt{-3}}(1,1,1,1,1,1)$; \item $\frac1{\sqrt{-3}}(1,1,1,\om,\om,\om)$ to
$-\frac1{\sqrt{-3}}(\om,\om,\om,1,1,1)$; \item $\frac1{\sqrt{-3}}(1,1,\om,\om,\om^2,\om^2)$ to itself.
\end{itemize}
After permutation, negation, and multiplication by $\om$, this takes care of virtually all cases.

Let \begin{equation}\label{fm}p_m(x_1,\ldots,x_6)=\sum_{(v_1,\ldots,v_6)}(v_1x_1+\cdots+v_6x_6)^m,\end{equation} where the sum is taken over the
756 vectors described above. Then $p_m$ is invariant under $G_{34}$ for every positive integer $m$. It is proved in \cite[Thm.10]{CS} that the
ring of complex invariant polynomials is given by
\begin{equation}\label{polyalg} \C[x_1,\ldots,x_6]^{G_{34}}=\C[p_6,p_{12},p_{18},p_{24},p_{30},p_{42}].
\end{equation}

In \cite{CS}, several other  lattices isomorphic to the above one are described, any of which can be used to give a different set of vectors $v$
and invariant polynomials $p_m$, still satisfying (\ref{polyalg}). The one that we have selected seems to give the simplest polynomials; in
particular, the only ones with integer coefficients.

We have $p_{6k}=S_1+S_2$, where $S_1=\sum_{i\ne j}\sum_{a,b=0}^2(\om^a x_i-\om^bx_j)^{6k}$, with $1\le i,j\le6$, and
$$S_2=\frac2{(-3)^{3k}}\sum_{a_i=0}^2(\om^{a_1}x_1+\cdots+\om^{a_5}x_5+\om^{-a_1-\cdots-a_5}x_6)^{6k}.$$
The coefficient of 2 on $S_2$ is due to the $\pm1$. Note that the sum for $S_1$ has $6\cdot5\cdot3^2$ terms, while that for $S_2$ has $3^5$
terms. Next note that if a term $T^{6k}$ occurs in either sum, then so does $(\om T)^{6k}$ and $(\om^2 T)^{6k}$, and all are equal. Thus we
obtain $S_1=3\sum_{i\ne j}\sum_{b=0}^2(x_i-\om^bx_j)^{6k}$ and
$$S_2=3\frac2{(-3)^{3k}}\sum_{a_2,\ldots,a_5=0}^2(x_1+\om^{a_2}x_2+\cdots+\om^{a_5}x_5+\om^{-a_2-\cdots-a_5}x_6)^{6k}.$$
We simplify $S_1$ further as
\begin{eqnarray*}S_1&=&3\sum_{\ell=0}^{6k}(-1)^\ell\tbinom{6k}\ell\sum_{i\ne j}x_i^\ell x_j^{6k-\ell}\sum_{b=0}^2\om^{b\ell}\\
&=&9\sum_{s=0}^{2k}(-1)^s\tbinom{6k}{3s}\sum_{i\ne j}x_i^{3s}x_j^{6k-3s}\\
&=&18(5m_{6k}+\sum_{s=1}^k(-1)^s\tbinom{6k}{3s}m_{(6k-3s,3s)}).\end{eqnarray*} At the first step, we have used that $\sum_{b=0}^2\om^{b\ell}$
equals 0 if $\ell\not\equiv0$ mod 3, and equals 3 if $\ell\equiv0$ mod 3. At the second step, we have noted that $\sum_{i\ne
j}x_i^{3s}x_j^{6k-3s}$ equals $m_{(6k-3s,3s)}$ if $s\not\in\{0,k,2k\}$, it equals $2m_{(3k,3k)}$ if $s=k$, and equals $5m_{(6k)}$ if $s=0$ or
$2k$.

The sum $S_2$ becomes
\begin{eqnarray*}S_2&=&\frac6{(-3)^{3k}}\sum_{\ebar}(\ebar)\sum_{a_2=0}^2(\om^{e_2-e_6})^{a_2}\cdots
\sum_{a_5=0}^2(\om^{e_5-e_6})^{a_5}x_1^{e_1}\cdots x_6^{e_6}\\
&=&\frac6{(-27)^k}\sum_{e_1\equiv\cdots\equiv e_6\ (3)}(\ebar)3^4x_1^{e_1}\cdots x_6^{e_6}.\end{eqnarray*}

Then $(-27)^k(S_1+S_2)/486$ equals the expression which we have listed for $f_{6k}$ in the statement of the theorem. We have chosen to work
with this rather than $p_{6k}$ itself for numerical simplicity. It is important that the omitted coefficient is not a multiple of 7.

For good measure, we show that (\ref{fm}) is 0 if $m\not\equiv0\ (6)$.
If $m\not\equiv0$ mod 3, then replacing terms $T^m$ by $(\om T)^m$ leaves the sums like $S_1$ and $S_2$ for (\ref{polyalg})
unchanged while, from a different perspective, it multiplies
them by $\om^m$. Thus the sums are 0. If $m\equiv 3$ mod 6, the term in $S_1$ corresponding to $\sum x_i^{3s}x_j^{m-3s}$ occurs with opposite
sign to that corresponding to $\sum x_i^{m-3s}x_j^s$, and so $S_1=0$. For $S_2$, the $(\pm1)^m$ will cause pairs of terms to cancel.
\end{pf*}

\begin{rmk}{\rm The only other place known to the author where formulas other than (\ref{fm}) for these polynomials exist is \cite{Ben}, where they
occupy 190 pages of dense text when printed.}\end{rmk}

As pointed out
by Kasper Andersen, $f_{42}-(f_6)^7$ is divisible by 7. This is  easily seen by expanding $(f_6)^7=(\sum(v_1x_1+\cdots+v_6x_6)^6)^7$
by the multinomial theorem. The need for this became apparent to Andersen, as the author had thought that the invariant ring of $G_{34}$
over $\Zh$ was $\Zh[f_6,\ldots,f_{42}]$, and this would have led to an impossible conclusion for the Adams operations in $QK^1(X_{34};\Zh)$.

Let $h_{42}=\frac17(f_{42}-(f_6)^7)$. Then we have the following result, for which we are grateful to Andersen.
\begin{thm}\label{Kgen} The invariant ring of $G_{34}$ over $\Zh$ is given by
$$\Zh[x_1,\ldots,x_6]^{G_{34}}=\Zh[f_6,f_{12},f_{18},f_{24},f_{30},h_{42}].$$
\end{thm}
\begin{pf} A {\tt Magma} program written and run by Andersen showed that each of these asserted generators is indecomposable
over $\Z/7$. (This is what failed when $f_{42}$ was used; it equals $(f_6)^7$ over $\Z/7$.) Thus the result follows from (\ref{polyalg}).\end{pf}

Since $f_{36}$ is invariant under $G_{34}$, it follows from (\ref{polyalg}) that it can be decomposed over $\C$ in terms of $f_6$, $f_{12}$,
$f_{18}$, $f_{24}$, and $f_{30}$. The nature of the coefficients in this decomposition was not so clear. It turned out that all coefficients
were rational numbers which are 7-adic units. We make this precise in
\begin{thm}\label{decompthm} $f_{36}$ can be decomposed as
\begin{eqnarray*}&&q_1f_6f_{30}+q_2f_{12}f_{24}+q_3f_{18}^2+q_4f_6^2f_{24}+q_5f_6f_{12}f_{18}\\
&&+q_6f_{12}^3+q_7f_6^3f_{18} +q_8f_6^2f_{12}^2+q_9f_6^4f_{12}+q_{10}f_6^6\end{eqnarray*}  with
\begin{eqnarray*}
q_1&=&944610925401/15161583716\\
q_2&=&733671261/19519520\\
q_3&=&243068633/9781739\\
q_4&=& -133840666859131062549/73986709144034080\\
q_5&=&-1758887990521258018071215403/629320589839873719708800\\
q_6&=&-1602221942044323/4879880000000\\
q_7&=&4011206338081535787030788541/114421925425431585401600\\
q_8&=&701461342458322269763709951654931/15733014745996842992720000000\\
q_9&=&-11844219519446025955021712628669/22348032309654606523750000\\
q_{10}&=&26589469730264682368719198549833/22348032309654606523750000
\end{eqnarray*}
Each of these coefficients $q_i$ is a 7-adic unit; i.e. no numerator or denominator is divisible by 7.
\end{thm}
\begin{pf} The ten products, $f_6f_{30},\ldots,f_6^6$, listed above are the only ones possible.
We express each of these products as a combination of monomial symmetric polynomials $m_\ebar$. We use {\tt Magma} to do this. The length of
$m_{(e_1,\ldots,e_r)}$ is defined to be $r$. We only kept track of components of the products of length $\le4$.  This meant that we only had to
include components of length $\le4$ of the various $f_{6k}$ being multiplied.

There were 34 $m_\ebar$'s of length $\le4$. These correspond to the partitions of 36 into multiples of 3. (Note that monomials with subscripts $\equiv$ 1 or
2 mod 3 only occur for us if the length is 6. Not having to deal with them simplifies our work considerably.) Indeed, there was one of length 1,
six of length 2, twelve of length 3, and fifteen of length 4. {\tt Magma} expressed each monomial such as $f_6f_{30}$ or $f_6^6$ as an integer
combination of these, plus monomials of greater length.  We just ignored in the output all those of greater length.  The coefficients in these
expressions were typically 12 to 15 digits. We also wrote $f_{36}$ as a combination of monomial symmetric polynomials of length $\le 4$,
ignoring the longer ones. This did not require any fancy software, just the multinomial coefficients from Theorem \ref{polys}.

Now we have a linear system of 34 linear  equations with integer coefficients in 10 unknowns. The unknowns are the coefficients $q_i$ in the
equation at the beginning of
\ref{decompthm}, and the
equations are the component monomials of length $\le 4$. Miraculously, there was a unique rational solution, as given in the statement of this theorem.

If it were not for the fact that the Conway-Sloane theorem \ref{polyalg} guarantees that there must be a solution when all monomial components
(of length $\le6$) are considered, then we would have to consider them all, but the fact that we got a unique solution looking at only the
monomial components of length $\le 4$ implies that this solution will continue to hold in the other unexamined components.
\end{pf}

Next we wish to modify the generators in \ref{Kgen} to obtain generators of $QK^1(X_{34};\Zh)$.
Similarly to \ref{29Kpolys}, we let $\ell_0(x)=\ln(1+x)$, and
\begin{equation}\label{Fidef}F_i=F_i(x_1,\ldots,x_6)=f_i(\ell_0(x_1),\ldots,f_6(\ell_0(x_6))).\end{equation}

A major calculation is required to modify the classes $F_i$ so that their coefficients are in $\Zh$; i.e. they do not have 7's in the
denominators. As observed after (\ref{AHSS}), it will be enough to accomplish this through grading 42 (with grading of $x_i$ considered to be 1).
\begin{thm} \label{thru42}The following expressions are 7-integral through grading 42:
\begin{itemize}
\item $F_{30}+\frac57F_{36}+\frac{22}{7^2}F_{42}$; \item $F_{24}+\frac47F_{30}+\frac{45}{7^2}F_{36}+\frac{104}{7^3}F_{42}$; \item
$F_{18}+\frac37F_{24}+\frac{20}{7^2}F_{30}+\frac{157}{7^3}F_{36}+\frac{526}{7^4}F_{42}$; \item
$F_{12}+\frac27F_{18}+\frac{45}{7^2}F_{24}+\frac{109}{7^3}F_{30}+\frac{1391}{7^4}F_{36}+\frac{6201}{7^5}F_{42}$; \item
$F_6+\frac17F_{12}+\frac{22}{7^2}F_{18}+\frac{204}{7^3}F_{24}+\frac{1107}{7^4}F_{30}+\frac{9682}{7^5}F_{36}+\frac{100682}{7^6}F_{42}$.
\end{itemize}
\end{thm}
It was very surprising that just linear terms were needed here. Decomposable terms were certainly expected. The analogue for $G_{29}$ in
\ref{29Kpolys} involved many decomposables. It would be interesting to know why Theorem \ref{thru42} works with just linear terms;
presumably this pattern will continue into higher gradings.
\begin{pf*}{Proof of Theorem \ref{thru42}}
 Similarly to the proof of \ref{29Kpolys}, we define
$$\Ft_i=\Ft_i(x_1,\ldots,x_6)=f_i(\ell_p(x_1),\ldots,\ell_p(x_6)),$$
and observe that a polynomial in the $\Ft_i$'s is 7-integral if and only if the same polynomial in the $F_i$'s is.

Next note that in the range of concern for Theorem \ref{thru42} $\ell_7(x)=x+x^7/7$. If we define
$$h_i=h_i(x_1,\ldots,x_6)=f_i(x_1+x_1^7,\ldots,x_6+x_6^7),$$
then \ref{thru42} is clearly equivalent to
\begin{stmt}\label{Geqs} For $t\ge1$ and grading $\le 42$,
\begin{itemize}
\item $h_{30}+5h_{36}+{22}h_{42}\equiv0$ mod $7^{t}$ in grading $30+6t$; \item $h_{24}+4h_{30}+{45}h_{36}+{104}h_{42}\equiv0$ mod $7^{t}$ in
grading $24+6t$; \item $h_{18}+3h_{24}+{20}h_{30}+{157}h_{36}+{526}h_{42}\equiv0$ mod $7^{t}$ in grading $18+6t$; \item
$h_{12}+2h_{18}+{45}h_{24}+{109}h_{30}+{1391}h_{36}+{6201}h_{42}\equiv0$ mod $7^{t}$ in grading $12+6t$; \item
$h_6+h_{12}+{22}h_{18}+{204}h_{24}+{1107}h_{30}+{9682}h_{36}+{100682}h_{42}\equiv0$ mod $7^{t}$ in grading $6+6t$.
\end{itemize}
\end{stmt}

We use {\tt Maple} to verify \ref{Geqs}. Our $f_i$'s are given in Theorem \ref{polys} in terms of $m_\ebar$'s. To evaluate
$m_\ebar(x_1+x_1^7,\ldots,x_6+x_6^7)$, the following result keeps the calculation manageable (e.g. it does not involve a sum over all
permutations). Partitions can be written either in increasing order or decreasing order; we use increasing. If $(a_1,\ldots,a_r)$ is an
$r$-tuple of positive integers, let $s(a_1,\ldots,a_r)$ denote the sorted form of the tuple; i.e. the rearranged version of the tuple so as to
be in increasing order. For example, $s(4,2,3,2)=(2,2,3,4)$.
\begin{prop}\label{mG} The component of $m_{(e_1,\ldots,e_r)}(x_1+x_1^7,\ldots,x_6+x_6^7)$ in grading $\sum e_i+6t$ is
$$\sum_\jbar \frac{P(e_1+6j_1,\ldots,e_r+6j_r)}{P(e_1,\ldots,e_r)}\binom{e_1}{j_1}\cdots\binom{e_r}{j_r}m_{s(e_1+6j_1,\ldots,e_r+6j_r)},$$
where $\jbar=(j_1,\ldots,j_r)$ ranges over all $r$-tuples of nonnegative integers summing to $t$, and $P(a_1,\ldots,a_r)$ is the product of the
factorials of repetend sizes.\end{prop} For example, $P(4,2,3,3)=2!$ because there are two 3's, $P(3,1,3,3,1,2)=3!2!$, and $P(3,4,2,1)=1$.

 \begin{expl}\label{expl} We consider as a typical example, the component of $$m_{(3,3,9,15)}(x_1+x_1^7,\ldots,x_6+x_6^7)$$ in grading 42.
 Table \ref{expltbl} lists the possible values of $\jbar$ and the contribution to the sum. The final answer is the sum of
 everything in the right hand column.
 \end{expl}

 \begin{center}
 \begin{minipage}{6.5in}
\begin{tab}{Terms for Example \ref{expl}}
\label{expltbl}
\begin{center}
\begin{tabular}{c|c}
$\jbar$&\text{term}\\[2pt]
\hline
$(2,0,0,0)$&$\binom32m_{3,9,15,15}$\\[3pt]
$(0,2,0,0)$&$\binom32m_{3,9,15,15}$\\[3pt]
$(0,0,2,0)$&$\binom92m_{3,3,15,21}$\\[3pt]
$(0,0,0,2)$&$\binom{15}2m_{3,3,9,27}$\\[3pt]
$(1,1,0,0)$&$3\cdot3\cdot3m_{9,9,9,15}$\\[2pt]
$(1,0,1,0)$&$3\cdot9m_{3,9,15,15}$\\[2pt]
$(1,0,0,1)$&$3\cdot15m_{3,9,9,21}$\\[2pt]
$(0,1,1,0)$&$3\cdot9m_{3,9,15,15}$\\[2pt]
$(0,1,0,1)$&$3\cdot15m_{3,9,9,21}$\\[2pt]
$(0,0,1,1)$&$9\cdot15m_{3,3,15,21}$
\end{tabular}
\end{center}
\end{tab}
\end{minipage}
\end{center}

\begin{pf*}{Proof of Proposition \ref{mG}} $m_{(e_1,\ldots,e_r)}(x_1+x_1^7,\ldots,x_6+x_6^7)$ is related to
\begin{equation}\label{sumsig}
\sum_\sigma(x_{\sigma(1)}^{e_1}+\tbinom{e_1}1x_{\sigma(1)}^{e_1+6}+\tbinom{e_1}2x_{\sigma(1)}^{e_1+12}+\cdots)\cdots(x_{\sigma(r)}^{e_r}+\tbinom
{e_r}1x_{\sigma(r)}^{e_r+6}+\cdots)\end{equation} summed over all permutations $\sigma$ in $\Sigma_r$. If $t$ values of $e_i$ are equal, then
(\ref{sumsig}) will give $t!$ times the correct answer. That is the reason that we divide by $P(\ebar)$. If $(e_1+6j_1,\ldots,e_r+6j_r)$
contains $s$ equal numbers, then the associated $m$ will be obtained from each of $s!$ permutations, which is the reason that
$P(e_1+6j_1,\ldots,e_r+6j_r)$ appears in the numerator.\end{pf*}

At first, mimicking \ref{29Kpolys}, we were allowing for products of $h$'s in addition to the linear terms which appear in \ref{Geqs}, but it
was turning out that what was needed to satisfy the congruences was just the linear term. If just a linear term was going to work, the
coefficients could be obtained by just looking at monomials of length 1. They were computed by {\tt Maple}, using that, by \ref{polys} and \ref{mG},
the coefficient of $m_{(6k+6t)}$ in $h_{6k}$ is
$(1+(-1)^k27^{k-1}\cdot5)\binom{6k}t$. Write the $k$th expression from the bottom of \ref{Geqs} as $\sum_{j\ge 0}a_{j,k}h_{6k+6j}$.
We require that the coefficient of $m_{(6k+6t)}$ in $\sum_{j=0}^ta_{j,k}h_{6k+6j}$ is 0 mod $7^t$. But this coefficient equals
$$\sum_{j=0}^ta_{j,k}\tbinom{6k+6j}{t-j}(1+(-1)^{k+j}27^{k+j-1}\cdot5).$$
We solve iteratively for $a_{j,k}$, starting with $a_{0,k}=1$, and obtain the values in \ref{Geqs}.
 Note that it first gives $a_{1,1}\equiv 1$ mod 7. If we had chosen
a value such as 8 or $-6$ instead of 1, then the value of $a_{2,1}$ would be different than 22. So these numbers $a_{j,k}$ are not uniquely
determined.  These different choices just amount to choosing a different basis for $QK^1(X_{34};\Zh)$.

Verifying Statement \ref{Geqs} required running many {\tt Maple} programs. For each line of \ref{Geqs}, a verification had to be made for each
relevant $t$-value, from two $t$-values for the first line down to six $t$-values for the last line. Moreover, for each of these pairs
(line number, $t$-value), it was convenient to use a separate program for monomials of each length 2, 3, 4, and 5, and then, for monomials of length 6,
it was
done separately for those with subscripts congruent to 0, 1, or 2 mod 3.  Thus altogether $(2+3+4+5+6)(4+3)=140$ {\tt Maple} programs were run. The
programs had enough similarity that one could be morphed into another quite easily, and a more skillful programmer could incorporate them all into the
same program.

Note that expanding from $f_j$ to $h_j$ does not change the number of components in monomials, nor does it change the mod 3 value of the
sum of the subscripts (i.e. exponents) in the monomials. This is simpler than the situation in the proof of \ref{decompthm}.
The algorithm is quite easy.
For each combination of $h$'s in \ref{Geqs}, replace each $h_{6j}$ by the combination of $m_{\ebar}$'s in $f_{6j}$ in \ref{polys}, but expanded using
\ref{mG}.
\end{pf*}

To obtain the Adams operations in $QK^1(X_{34};\Zh)$, we argue similarly to the paragraph which precedes Theorem \ref{psik}. First note that $F_{36}$ decomposes in terms of $F_{6i}$'s exactly as
does $f_{36}$ in terms of $f_{6i}$'s in \ref{decompthm}. We can\footnote{But we need not bother to do so explicitly.} modify by decomposables
in dimensions greater than 42
to obtain  7-integral classes $G_6$, $G_{12}$, $G_{18}$, $G_{24}$, and $G_{30}$ which
agree with the classes of \ref{thru42} (with $F_{36}$ decomposed) through dimension 42. There is also a 7-integral class $G_{42}$ which agrees
with $\frac17(F_{42}-(F_6)^7)$ in dimension 42. These generate $K^*(BX_{34};\Zh)$ as a power series algebra. As in the preamble to \ref{psik},
then $z_i:=B^{-1}e^*(G_{i+1})$ for
$i=5$, 11, 17, 23, 29, and 41 form a basis for $QK^1(X_{34};\Zh)$, and $e^*$ annihilates decomposables.

Similarly to the situation for $(X_{29})_5$ in the proof of \ref{psik}, if we let
$$P=\begin{pmatrix}1&0&0&0&0&0\\[4pt]
 \frac17&1&0&0&0&0\\[4pt]
  \frac{22}{49}&\frac27&1&0&0&0\\[4pt]
\frac{204}{343}&\frac{45}{49}&\frac37&1&0&0\\[4pt]
\frac{1107}{2401}&\frac{109}{343}&\frac{20}{49}&\frac47&1&0\\[4pt]
\frac{16647}{16807}&\frac{1399}{2401}&\frac{183}{343}&\frac{6}{49}&\frac17&1\end{pmatrix},$$
then the matrix of $\psi^k$ on the basis $\{z_5,z_{11},z_{17},z_{23},z_{29},z_{41}\}$ is
$$P^{-1}\text{diag}(k^5,k^{11},k^{17},k^{23},k^{29},k^{41})P.$$ The entries in the last row of $P$ are 7 times the coefficients of $F_{42}$ in
\ref{thru42} reduced mod 1. Those coefficients were multiplied by 7 because $z_{41}$ is related to $\frac17F_{42}$ rather than to $F_{42}$.

Using this, we compute the $v_1$-periodic homotopy groups, similarly to \ref{v*G29}. Note the remarkable similarity with that result. Here, of
course, $\nu(-)$ denotes the exponent of 7 in an integer.
\begin{thm}\label{Adams} The groups $\vp_*(X_{34})_{(7)}$ are given by
$$\vp_{2t-1}(X_{34})\approx\vp_{2t}(X_{34})\approx\begin{cases}0&t\not\equiv5\ (6)\\
\Z/7^5&t\equiv 5,35\ (42)\\
\Z/7^{\min(12,5+\nu(t-11-12\cdot 7^6))}&t\equiv 11\ (42)\\
\Z/7^{\min(18,5+\nu(t-17-18\cdot 7^{12}))}&t\equiv 17\ (42)\\
\Z/7^{\min(24,5+\nu(t-23-18\cdot 7^{18}))}&t\equiv 23\ (42)\\
\Z/7^{\min(30,5+\nu(t-29-12\cdot 7^{24}))}&t\equiv 29\ (42)\\
\Z/7^{\min(42,5+\nu(t-41-24\cdot 7^{36}))}&t\equiv 41\ (42).\end{cases}$$
\end{thm}

\begin{pf} The group $\vp_{2t}(X)_{(7)}$ is presented by $\begin{pmatrix}(\psi^7)^T\\(\psi^3)^T-3^tI\end{pmatrix}$,
since 3 generates $\Z/49^\times$. We let $x=3^t$ and form this matrix analogously to (\ref{psimat}). Five times we can pivot on units, removing
their rows and columns, leaving a column matrix with 7 polynomials in $x$. The 7-exponent of $\vp_{2t}(X)_{(7)}$ is the smallest of that of
these polynomials (with $x=3^t$). This will be 0 unless $x\equiv5$ mod 7, which is equivalent to $t\equiv5$ mod 6. We find that two of these polynomials will
always yield, between them, the smallest exponent. Similarly to (\ref{peq}) and Table \ref{exptbl}, we write these polynomials as $p_i(3^m+y)$ for
carefully-chosen values of $m$. Much preliminary work is required to discover these values of $m$. Ignoring unit coefficients and ignoring
higher-power terms whose coefficients will be sufficiently divisible that they will not affect the divisibility, these polynomials will be as in
Table \ref{polystbl}.

\begin{tab}{Certain  $p_i(3^m+y)$, (linear part only)}
\label{polystbl}
\begin{center}
\begin{tabular}{c|cc}
$m$&$p_1$&$p_2$\\
\hline
$5,\ 35$&$7^5+7^4y$&$7^{\ge12}+7^{11}y$\\[2pt]
$11+12\cdot7^6$&$7^{12}+7^4y$&$7^{12}+7^{11}y$\\[2pt]
$17+18\cdot7^{12}$&$7^{18}+7^4y$&$7^{18}+7^{11}y$\\[2pt]
$23+18\cdot7^{18}$&$7^{25}+7^4y$&$7^{24}+7^{11}y$\\[2pt]
$29+12\cdot7^{24}$&$7^{31}+7^4y$&$7^{30}+7^{11}y$\\[2pt]
$41+24\cdot7^{36}$&$7^{43}+7^4y$&$7^{42}+7^{11}y$
\end{tabular}
\end{center}
\end{tab}

The claim of the theorem follows from Table \ref{polystbl} by the same argument as was used in the proof
of \ref{v*G29}. For $t$ in the specified congruence, if $3^t=3^m+y$, then $\nu(y)=\nu(t-m)+1\ge2$, similarly to (\ref{2025}).
For example, if $t\equiv11$ mod 42, and $3^t=3^{11+12\cdot7^6}+y$, then $\nu(y)=\nu(t-11-12\cdot7^6)+1$.
Thus $\min(\nu(p_1(3^t)),\nu(p_2(3^t)))$ will be determined by the $7^5$ in $p_1$ if $t\equiv5,35\ (42)$,
while in the other cases, it is determined by the $7^4y$ in $p_1$ or the constant term in $p_2$.

The groups $\vp_{2t-1}(X_{34})$ are cyclic by an argument  similar to the
one described at the end of the proof of \ref{v*G29}, and have the same order as $\vp_{2t}(X_{34})$
for the standard reason described there.
\end{pf}

Similarly to the discussion preceding Theorem \ref{29conj}, one of the factors in the product decomposition of
$SU(42)_7$ given in \cite{MNT} is an $H$-space $B_5^7(7)$ whose $\F_7$-cohomology is an exterior algebra on classes of
grading 11, 23, 35, 47, 59, 71, and 83, and which is built from spheres of these dimensions by fibrations.
Using \cite{Ya}, we can obtain a degree-1 map $B_5^7(7)\to S^{71}$. Let $B_7:=B(11,23,35,47,59,83)$ denote its fiber.
The following result was conjectured by the author and proved by John Harper. Its proof will be described in the last line of the paper.
\begin{thm}\label{34conj} (Harper) There is a homotopy equivalence $(X_{34})_7\simeq B_7$.\end{thm}

\section{Proofs provided by John Harper}\label{Harsec}
In this section, we provide proofs of Theorems \ref{29conj}, \ref{X31E8}, and \ref{34conj}, which were explained to
the author by John Harper. We begin with the following strengthening of Theorem \ref{CHZthm}.
\begin{thm}\label{rankp-1} Theorem \ref{CHZthm} is true for $r\le p-1$.\end{thm}
\begin{pf} Although the proof presented in \cite{CHZ} works for $r=p-1$ just as it does for $r<p-1$, we take this opportunity
to explain some aspects of it more thoroughly. For $\l\ne0\in\Z/p$, a $\pi_\l$-space (resp. $Q_\l$-space) is one which admits a self-map inducing multiplication
by $\l$ in $\pi_*(-)\otimes\Z/p$ (resp. $QH^*(-;\Z/p)$). A $\pi_\l$-map (resp. $Q_\l$-map) is a map between spaces of the indicated type
which commutes up to homotopy with the self-maps.

The first part of the proof, extending \cite[Thm 1.3]{CHZ}, involves showing that if $X$ is an $H$-space of rank $r\le p-1$ and with torsion-free
homology, then there is a $\pi_\l$-map $X\to S^n$ whose homotopy fiber $Y$ is an $H$-space of rank $r-1$ with torsion-free homology.
The proof of this is exactly as in \cite[p.358]{CHZ}, noting that \cite[4.2.2]{Zab} is valid (and stated) for $r=p-1$ as well as for
$r<p-1$. That $Y$ is an $H$-space follows from \cite[Thm 1.1]{CHZ} since its rank is less than $p-1$. This construction can be iterated
to yield (\ref{fibrs}).

We use Lemma \ref{XY} to deduce the part of \ref{rankp-1} which says that the homotopy type of $X$ is determined by certain homotopy classes.
Suppose we have (\ref{fibrs}) and a primed version. Suppose we have shown that
$X_i$ and $X_i'$ are $p$-equivalent $H$-spaces and that the elements  $\a\in\pi_{n_{i+1}-1}(X_i)$ and $\a'\in\pi_{n_{i+1}-1}(X'_i)$ correspond under
this equivalence. Then, in the notation of \ref{XY}, there are equivalences as $\pi_\l$-spaces
$$X_{i+1}\simeq (X_i)^\a\simeq (X_i')^{\a'}\simeq X_{i+1}'.$$
If $i+1<p-1$, then by \cite[Thm 1.1]{CHZ}, there are $p$-equivalent $H$-space structures on $X_{i+1}$ and $X'_{i+1}$, extending the induction.
Note that for $X_{p-1}$ and $X'_{p-1}$ we do not assert an equivalence as $H$-spaces, only as $\pi_\l$-spaces.
\end{pf}
\begin{lem}\label{XY} Suppose $Y\to X@>g>> S^n$ is a fibration, with $X$ a $\pi_\l$-space and $Y$ an $H$-space with rank$(Y)<p-1$.
Let $\a\in\pi_{n-1}(Y)$ denote $\partial(\iota_n)$ in the homotopy sequence of the fibration.
As in \cite[p.351]{CHZ}, let $Y^\a=D^n_+\times Y\cup_cD^n_-\times Y$, where $c(x,y)=(x,\a(x)y)$, using the $H$-space multiplication of $Y$.
Then there is a homotopy equivalence $X\simeq Y^\a$ which is a $\pi_\l$-map.\end{lem}
\begin{pf} We expand on some aspects of the proof given in \cite[p.358]{CHZ} and correct several confusing typos.
As explained there, $Y^\a$ is a $\pi_\l$-space
and admits an inclusion $Y\cup_\a e^n\hookrightarrow Y^\a$.

On the other hand, the given fibration is fiber homotopy equivalent to $Y\to X'\to S^n$, with
$X'=D^n_+\times Y\cup_\g D^n_-\times Y$, where $\g:S^{n-1}\times Y\to S^{n-1}\times Y$ is the clutching function defined
by the homotopy equivalences $D^n_\pm\times Y\to g^{-1}(D^n_\pm)$. We will work with $X'$, but in the end may replace it with $X$.
Note that $X'$ inherits the $\pi_\l$-structure of $X$.

Since $\g|S^{n-1}\times\{*\}$ represents $\a$, the relative $n$-skeleta $(Y^\a,Y)^n$ and $(X',Y)^n$ both equal $Y\cup_\a e^n$.
Thus the $n$th stages of the Postnikov systems of $Y^\a$ and $X'$ are built from the same Eilenberg-MacLane spaces and same
$k$-invariants, and so there exists a $\pi_\l$-map $(Y^\a)_n\to (X')_n$ of these $n$th stages, inducing a cohomology isomorphism
in dimension $\le n$. We wish to show that this extends to a $\pi_\l$-map $Y^\a\to X'$, which will then be a homotopy equivalence
by the Five Lemma applied to the homotopy sequences of the fibrations $Y\to Y^\a\to S^n$ and $Y\to X'\to S^n$, and Whitehead's Theorem.

This requires that we consider the  primitive Postnikov systems (PPS) of the map $X'\to (X')_n$.
The {\it primitive Postnikov system}, as described, for example, in \cite[p.426ff]{Kane}, applies to a rational equivalence, and
gives a Postnikov tower in which all fibers are mod $p$ Eilenberg-MacLane spaces. Our map $X'\to (X')_n$
is a rational equivalence because the only infinite homotopy groups of the spheres that build  $X'$ are present in
 $(X')_n$. By \cite[Lemma 2.8]{CHZ}, the maps in the PPS are $\pi_\l$-maps. We will also use that, by \cite[Lemma 2.4]{CHZ},
 all our $\pi_\l$-maps are also $Q_\l$-maps.

We will show that the composite $Y^\a\to (Y^\a)_n\to (X')_n$ lifts through the PPS of $X'\to (X')_n$.
Assume there is a lifting to a $\pi_\l$-map $Y^\a\to E_s$, where $E_s$ is some stage in the PPS.
Since the $k$-invariants in the PPS are $\l$-eigenvectors by \cite[Lemma 2.8]{CHZ}, so are their images in $H^*(Y^\a;\F_p)$.
A basis for $H^*(Y^\a;\F_p)$ consists of products of no more than $p-1$ generators, each of which is a $\l$-eigenvector
for the $Q_\l$-map of $Y^\a$. An element in this basis which is not one of the generators is a $\l^i$-eigenvector
for some $2\le i\le p-1$. Since, for such $i$, $\l^i\not\equiv \l$ mod $p$, we deduce that a nonzero image of a $k$-invariant can only
equal one of the generators\footnote{or a linear combination of generators in the same dimension} of $H^*(Y^\a;\F_p)$.
However, the $k$-invariants are in dimension greater than that of the
generators of $H^*(Y^\a;\F_p)$. We conclude that the image of the $k$-invariants must equal 0, and so the map lifts
to $Y^\a\to E_{s+1}$. By \cite[2.7]{CHZ}, this lifting may be chosen to be a $\pi_\l$-map. Since the PPS has only finitely many
stages through dim$(Y^\a)$, we obtain the desired lifting $Y^\a\to X'$.
\end{pf}

We now apply these general results to our specific situations.

\begin{pf*}{Proof of Theorem \ref{X31E8}.} Both spaces $(X_{31})_5$ and $X_0(E_8)$ are $H$-spaces with mod-5 cohomology an
exterior algebra on classes of dimension 15, 23, 39, and 47. By Theorem \ref{rankp-1}, there exist diagrams of fibrations
$$\begin{CD} S^{15}@>>> X_2@>>>X_3@>>>(X_{31})_5\\
@. @Vf_2VV @Vf_3VV @Vf_4VV\\
@.S^{23}@.S^{39}@.S^{47} \end{CD}$$
and
$$\begin{CD}S^{15}@>>>X_2'@>>>X_3'@>>>X_0(E_8)\\
@.@Vf'_2VV @Vf'_3VV @Vf'_4VV\\
@.S^{23}@. S^{39}@.S^{47},\end{CD}$$
and the homotopy types of the spaces are determined by relevant elements of homotopy groups.

Always localized at 5, we have $\pi_{22}(S^{15})\approx\Z/5$ generated by $\a_1$. The elements in $\pi_{22}(S^{15})$
which determine both $X_2$ and $X_2'$ are nonzero multiples of $\a_1$, and so by \ref{CHZthm} we obtain an equivalence $X_2'\simeq X_2$.
The claim about these homotopy classes is proven for $(X_{31})_5$ from the entry in position (2,1) in the Adams operation matrix
in \ref{31ops}; Adams' $e$-invariant says that $\a_1$ attaching maps are present iff the relevant Adams operation involves
$u(k^n-k^{n+p-1})/p$ with $u$ a unit. Similarly for $X_0(E_8)$, the presence of $\a_1$ can be deduced from the Adams operations
in $E_8$ as given in \cite[Prop 3.5]{Reprth}, although here the well-known Steenrod algebra action in $H^*(E_8)$ also implies
the attaching map. We now identify $X_2'$ with $X_2$ in our notation.

The homomorphism $\pi_{38}(X_2)\to \pi_{38}(S^{23})$ is a surjection $\Z/25\to\Z/5$. See, e.g., Diagram \ref{Bch}.
The $\frac15(k^{11}-k^{19})$ in position (3,2) of \ref{31ops} and in the formula for $\psi^k(x_{11})$ in \cite[Prop 3.5]{Reprth}
tell that in both of our sequences, the determining element of $\pi_{38}(X_2)$ is a generator, i.e., it maps to $\a_2\in\pi_{38}(S^{23})$,
and so, by Theorem \ref{CHZthm},
there is a homotopy equivalence $X_3\simeq X_3'$. This is probably the only place that any of our calculations with $(X_{31})_5$
are required in proving \ref{X31E8}. The $\a_1$ attaching maps can be seen by Steenrod operations, but $\a_2$ requires secondary operations,
which were used in \cite{Gon} to see the $\a_2$ attaching map in $X_0(E_8)$, or Adams operations in $K$-theory, as we have done.

We identify $X_3$ and $X_3'$. The homomorphism $\pi_{46}(X_3)\to \pi_{46}(S^{39})$ is a surjection $\Z/5^3\to \Z/5$. (The cyclicity of
$\pi_{46}(X_3)$ can be proved
 easily using \cite[Prop 5.5]{Reprth}.) As in the previous two paragraphs, Adams operations imply that our determining
element in both sequences is a generator, i.e., it maps to $\a_1\in\pi_{46}(S^{39})$, and so we obtain the asserted equivalence from Theorem \ref{rankp-1}. \end{pf*}

\begin{pf*}{Proof of Theorem \ref{29conj}.} We first show that the map $B_3^5\to S^{31}$ described prior to \ref{29conj} can be chosen to be
a $\pi_\l$-map. By \cite[Lemma 2.8]{CHZ}, the generating map $S^{31}\to Z(\Z,31)$ has a PPT in which all maps are $\pi_\l$-maps.
The map $B_3^5@>g>> K(\Z,31)$ is a $\pi_\l$-map, since $\pi_{31}(B_3^5)$ is cyclic. By \cite{Ya}, $g$ lifts to a map $B^5_3\to S^{31}$.
By the Lifting Theorem, \cite[2.7a]{CHZ}, this lifting, through each stage of the PPT, is a $\pi_\l$-map.

Thus, by \cite[Lemma 2.3]{CHZ}, the space which we call $B(7,15,23,39)$ in \ref{29conj} is a $\pi_\l$-space.
Let $f$ denote the composite $B(7,15,23,39)\to B_3^5\to SU(20)\to S^{39}$. Its fiber, $B_3^3=B(7,15,23)$, is one of the four
factors in $SU(12)_5$, and hence is an $H$-space.

On the other hand, Theorem \ref{rankp-1} yields fibrations of $H$-spaces, localized at 5,
\begin{equation}\label{Xfibrs}\begin{CD}S^7@>>> X_2@>>> X_3@>>> (X_{29})_5\\
@. @VVV @VVV @VVV\\
@.S^{15}@. S^{23}@.S^{39}\end{CD}\end{equation}
We can almost apply Theorem \ref{rankp-1} to obtain the desired equivalence of $(X_{29})_5$ and $B(7,15,23,39)$.
The only thing missing is that we do not know that $B(7,15,23,39)$ is an $H$-space.

Instead, we first use Theorem \ref{CHZthm} to show that $X_3$ (of (\ref{Xfibrs})) and $B_3^3$ are homotopy equivalent $H$-spaces.
Note that in the proof of \ref{rankp-1} we obtained equivalent $H$-structures as long as the rank was less than $p-1$.
In this case, the homotopy classes that must correspond are, in both cases, $\a_1\in\pi_{14}(S^7)_{(5)}\approx\Z/5$ and a generator
of $\pi_{22}(B(7,15))\approx\Z/25$, mapping to $\a_1\in\pi_{22}(S^{15})$. The $\a_1$'s in $B_3^3$ are well-known
by the action of $\P^1$, while in $X_3$ they may be deduced from entries in the matrix of \ref{psik} in positions $(2,1)$ and $(3,2)$.

For the fibrations $X_3\to (X_{29})_5\to S^{39}$ and $B_3^3\to B(7,15,23,39)\to S^{39}$, the relevant homotopy groups, always localized at 5,
$\pi_{38}(X_3)$ and $\pi_{38}(B_3^3)$, are both $\Z/5^3$, and the relevant classes are generators, which we will call $\a_2$ since they each map
to $\a_2\in\pi_{38}(S^{23})$.
This is seen for $X_3$ from the entry in position $(4,3)$ of the matrix of \ref{psik}, while for $B_3^3$ it can be obtained from the
Adams operation calculation in \cite[pp 667-668]{Pot}. Thus, by Lemma \ref{XY}, we have homotopy equivalences
$$(X_{29})_5\simeq (X_3)^{\a_2}\simeq (B_3^3)^{\a_2}\simeq B(7,15,23,39),$$
as desired.\end{pf*}

The proof of Theorem \ref{34conj} is entirely analogous, and is omitted.

\def\line{\rule{.6in}{.6pt}}

\end{document}